\documentclass[a4paper,10pt]{article}
\usepackage{amsfonts}
\usepackage{graphicx,amssymb,amsmath}
\usepackage{amscd}
\usepackage[english]{babel}
\usepackage[ansinew]{inputenc}
\usepackage{slashbox}
\usepackage[all]{xy}
\usepackage{latexsym}
\usepackage{amsmath}
\usepackage{stmaryrd}

\newcommand {\ds}{\displaystyle}

 \hsize \textwidth
 \advance \hsize by -\marginparwidth
\newcommand{\field}[1]{\mathbb{#1}}
\newcommand{\R}{\field{R}}

\newcommand{\N}{\field{N}}

\newtheorem{theo}{Theorem}[section]
\newtheorem{theor}{Theorem}
\newtheorem{obs}[theo]{Observations}

\newtheorem{defi}[theo]{Definition}
\newtheorem{lem}[theo]{Lemma}
\newtheorem{cor}[theo]{Corollary}
\newtheorem{prop}[theo]{Proposition}
\newtheorem{rem}[theo]{Remark}

\newtheorem{exe}[theo]{Example}
\newtheorem{exes}[theo]{Examples}

\def\dis{\ds}
\def\ap{\rightarrow}
\def\S{\Sigma}
\def\l{\lambda}
 \def\a{\alpha}
\def\b{\beta}
\def\g{\gamma}
\def\G{\Gamma}
\def\D{\triangle}
\def\t{\tau}
\def\d{\delta}
\def\o{\omega}

\def\l{\lambda}
\def\L{\Lambda}

\def\r{\rho}

\def\s{\sigma}

\def\so{\underline}

\def\O{\Omega}

\def\~{\tilde}
\def\dis{\displaystyle}
\input epsf

\title{Integrability of weak distributions on Banach manifolds}
\author{   F. Pelletier\footnote{Laboratoire de Math\'ematiques, Universit\'e de Savoie, Campus scientifique, 73376 Le Bourget du Lac, France}\footnote
{The author is grateful to
Prof Tilmann Wurzbacher  for long and helpful discussions about  integrability of distributions}
}
\date{}
\begin{document}
\maketitle

\begin{abstract}
This paper concerns the problem of integrability of non closed distributions on Banach manifolds. We introduce the notion of weak distribution and we  look for conditions under which these distributions admit weak integral submanifolds.  We give some applications to Banach Lie algebroids and Banach  Poisson manifold. The main  results of this paper  generalize the works presented in \cite{ChSt}, \cite {Nu} and \cite{Gl}. \end{abstract}

\bigskip

\noindent {\bf AMS Classification}: 58B10, 53C12, 53B50, 53B25, 47H10, 47A60, 46B20, 46B07, 37K25, 37C10. \\

\noindent {\bf Key words}: Banach manifold, weak Banach submanifold, weak distribution, integral manifold, involutive distribution, integrable distribution, invariance, Lie invariance, Banach Lie algebroid, Banach Poisson manifold.

%%%%%%%%%%%%%%%%%%%%%%%%
\section{Introduction}\label{intro}
%%%%%%%%%%%%%%%%%%%%%%%%%%
${}\;\;\;\;\;\;$ In  differential geometry, a {\it distribution} on a smooth manifold $M$ is an assignment ${\cal D}: x\mapsto {\cal D}_x\subset T_xM$ on $M$, where ${\cal D}_x$ is a subspace of $T_xM$. The distribution is {\it integrable} if, for any $x\in M$ there exists an immersed  submanifold $f:L\ap M$ such that $x$ belongs to  $f(L)$ and for any $z\in L$ we have $Tf(T_zL)={\cal D}_{f(z)}$. On the other hand, ${\cal D}$ is called {\it involutive} if, for any vector fields $X$ and $Y$  on $M$ which are tangent to $\cal D$, the Lie bracket $[X,Y]$ is also tangent to $\cal D$. The distribution is {\it invariant} if for any vector field $X$ tangent to ${\cal D}$, the flow $\phi^X_t$ leaves ${\cal D}$ invariant (see \ref{preliminaires et resultats}).\\

On  finite dimensional  manifod, when ${\cal D}$ is a subbundle of $TM$,
the classical Frobenius Theorem gives an equivalence between integrability and involutivity. In the other case, the distribution is "singular" and even under  assumptions of smoothness on ${\cal D}$, in general,  the involutivity is not a sufficient condition for integrability (we need some more additional  local conditions). These problems were clarified and resolved essentially in \cite{Su}, \cite{St} and \cite{Ba}.\\

In the context of Banach manifolds, the Frobenius Theorem is again true, for distributions which are  complemented subbundles  in  the tangent bundle. For singular distributions, some papers (\cite{ChSt}, \cite{Nu} for instance) show that, when  the distribution is closed and complemented (i.e. ${\cal D}_x$ is a complemented Banach subspace of $T_xM$), we have equivalence between integrability and invariance. Under  sufficient conditions about local involutivity  we also  get  a result of integrability. A more recent work (\cite{Gl}) proves analog results without the assumption that the distribution is complemented.\\

According to the notion of  "weak submanifolds" in a Banach manifold introduced in ((\cite{El},\cite{Pe}),  in this paper, we consider "weak distributions":    ${\cal D}_x$ can be not closed in $T_xM$ but  ${\cal D}_x$ is endowed with its own Banach structure, so that the inclusion ${\cal D}_x\ap T_xM$ is continuous. Such a category of distributions takes naturally place in the framework of Banach Lie  algebroids (morphisms from a Banach bundle over a Banach manifold into the tangent bundle of this manifold). Under conditions of "local triviality", our result can be seen like generalization as well of results of \cite{Su} and  \cite{St}, than of results of \cite {Nu} and \cite{Gl}.  Note that, our proofs take in account Remarks of Balan in \cite{Ba}, about results of \cite{Su} and \cite{St} (see Obervations \ref{balan})\\

The first section contains  the most important  definitions and properties about weak distributions. It contains  also the first result of equivalence between integrability  and invariance (Theorem 1), under local lower triviality assumption.  This last property is, in fact,  a generalization of the classical notion of "smoothness" for distributions (see Observations \ref{smoothG}). In the second section, we adapt  the arguments used in   \cite{ChSt}  to our context:  under condition of "Lie invariance", we give a generalization of their results  about  the  integrability of distributions (Theorem 2). In the second part, under the assumption of "upper triviality"  (which is a general context for anchored bundles), we give some  conditions of "local involutivity"  which gives rise to an integrability property (Theorem 4). In the last section, we give some applications of these results in the context of  Banach Lie algebroids  and Banach Poisson manifold ( cf \cite{OdRa1} and \cite{OdRa2})

.%%%%%%%%%%%%%%%%%%%%%%%%%%%%%%%%%%%%%%%%%%%%%%%%%
\section{Integrability and invariance}\label{intinv}
%%%%%%%%%%%%%%%%%%%%%%%%%%%%%%%%%%%%%%%%%%%%%%%%%
\subsection{Preliminaries and context}\label{preliminaires et resultats}
%%%%%%%%%%%%%%%%%%%%%%%%%%%%%%%%%%%%%%%%%%%%%%%%%%%
%\subsection{Preliminaries}\label{preliminaires}
%%%%%%%%%%%%%%%%
Let $M$ be  a smooth connected Banach manifold  modelled on a Banach space $E$. We denote by  $C^\infty(M)$ the ring of smooth  functions on  $M$  and by  ${\cal X}(M)$ the Lie algebra of smooth vector fields on $M$.
 A {\bf local vector field } $X$ on  $M$ is a smooth section of the tangent bundle  $TM$   defined on an open set of  $M$ (denoted by  Dom($X$)). Let be   ${\cal X}_L(M)$ the set of all local vector fields  on $M$. Such a vector field $X\in {\cal X}_L(M)$  has a flow  $\phi^X_t$ which is defined on a maximal  open set  $\O_X$ of  $M\times\R$.

 Using the terminology introduced in \cite{El} or \cite{Pe}, 
  a {\bf weak submanifold}  of $M$  is a pair $(N,f)$ where $N$ is  Banach manifold  modelled on a Banach space $F$ and $f:N\ap M$ a smooth map  such that: 
\begin{enumerate}
\item there exists an injective continuous  linear map $i:F\ap E$ between these two Banach spaces;
\item $f $ is injective and the tangent map $T_xf:T_xN\ap T_{f(x)}M$ is injective for all $x\in N$.\\
\end{enumerate}
\begin{rem}\label{weakN}${}$\\
 Given a weak submanifold $f:N\ap M$, on the subset $f(N)$ in $M$ we have two topologies:
 \begin{enumerate}
   \item the induced topology from $M$
   \item the topology for which $f$ is a homeomorphism from $N$ to $f(N)$.
 \end{enumerate}
 With this last topology, via $f$, we get a structure of Banach manifold modelled on $F$.
 Moreover, the inclusion from $f(N)$ into $M$ is continuous as map from the Banach manifold $f(N)$ to $M$.
 In particular, if  $U$ is  open in $M$, then, $f(N)\cap U$ is an open set for the topology of the Banach manifold  on $f(N)$.\\
 Note that in \cite{El} and \cite{Pe} the  definition of a "weak manifold"  imposes that these two topologies are identical. Our definition is somewhat different and is motivated by the notion of  "weak immersion" introduced in \cite {OdRa1} and \cite{OdRa2}  
 \end{rem}

In this work, a  {\bf weak distribution} on  $M$  is a map $ {\cal D}: x\ap {\cal D}_x$  which, for every $x\in M$, associates  a vector subspace  ${\cal D}_x$ in $T_xM$  (not necessarily closed) endowed with a norm $||\;||_x$ so that   $({\cal D}_x, ||\;||_x)$ is a Banach space (denoted by $\tilde{\cal D}_x$) and such  that the natural inclusion $i_x : \tilde{\cal D}_x \ap T_xM$ is continuous.

\begin{rem}\label{closed}${}$\\
When ${\cal D}_x$ is closed, via any chart, we get a norm on $T_xM$ which induces a Banach structure on ${\cal D}_x$. So if ${\cal D}_x$ is closed for all $x\in M$, the previous assumption on the Banach structure $\tilde{\cal D}_x$ is always satisfied, and we get the usual definition of a distribution on $M$ (compare with \cite{Gl}, \cite{ChSt}, \cite{Nu}). In this last situation we always endow $\tilde{\cal D}_x$ with this induced Banach structure and we say that $\cal D$ is {\bf closed}.
\end{rem}

\begin{exes}\label{ex}${}$
\begin{enumerate}
\item Let $l^p$ (resp. $l^\infty$) be the Banach space of real  sequences $(x_k)$ such that $\ds\sum_{k=1}^\infty |x_k|^p<\infty$ (resp. absolutely bounded) and denote by $I_p$ the natural inclusion of $l^1$ in $l^p$ , $p>1$ or $p=\infty$. On the Banach space $l^p$, $x\mapsto {\cal D}_x=x+I_p(l^1)$ is a weak  distribution which is not closed.
\item Let $E$ and $F$ be two Banach spaces and $T:F\ap E$ a continuous operator. Denote by $\hat{T}:F/\ker T\ap E$ the canonical quotient bijection associated to $T$  that is
\begin{eqnarray}\label{quotient}
\xymatrix {
    F  \ar[r]^{^q\;\;\; \;\;\;\;}\ar[d]_{_{T}} & {F/\ker{T}} \ar[ld]^{\hat{T}} \\
   T(F)}
 \end{eqnarray}
We can endow $T(F)$ with the structure of Banach space such that $\hat{T}$ is an isometry.  On $E$,   the assignment   $x\mapsto {\cal D}_x=x+T(F)$ is a weak distribution. This distribution is closed if and only if $T(F)$ is closed in $E$.
\item Let $L(F,E)$ be the set of continuous operators between the Banach spaces $F$ and
$E$. Given  a smooth map $\Psi: E\ap L(F,E)$, we denote by $\Psi_x$ the continuous operator associated to $x\in E$. As  in  (\ref{quotient}) denote by  $\hat{\Psi}_x$ the canonical bijection associated $\Psi_x$ and we endow ${\cal D}_x=\Psi_x(F)$ with the Banach structure such that $\hat{\Psi}_x$ is an isometry. Then, $x\ap {\cal D}_x$ is a weak distribution on $E$ which is closed if and only if ${\cal D}_x$ is closed for any $x\in E$.
\end{enumerate}
\end{exes}

\noindent A vector field $Z\in{\cal X}(M)$ is {\bf tangent} to $\cal D$, if for  all $x\in$ Dom($Z$), $Z(x)$ belongs to ${\cal D}_x$. The set of local vector fields tangent to $\cal D$  will be denoted by {\bf $\bf {\cal X}_{\cal D}$}.

We say that  $\cal D$ is {\bf is generated by a subset  }  ${\cal X}\subset {\cal X}_L(M)$ if, for every $x\in M$, the vector space ${\cal D}_x$ is the  linear hull  of the set $\{Y(x)\;,\;Y\in{\cal X}\;,\;x\in$ Dom$(Y)\}$.  \\

For a weak distribution $\cal D$ on $M$, we have the following classical  properties:\\

$\bullet$ an {\bf integral manifold} of $\cal D$  through $x$ is a weak  submanifold $f:N\ap M$ such that $f(\tilde{x})=x$  for some  $\tilde{x}\in N$  and  and $T_{\tilde{y}}f(T_{\tilde{y}}N)={\cal D}_{f(\tilde{y})}$ for all $\tilde{y}\in N$.\\

$\bullet$   $\cal D$ is called {\bf integrable} if for any $x\in M$ there exists an integral manifold $N$ of $\cal D$ through $x$.\\

$\bullet$ if  $\cal D$ is generated by ${\cal X}\subset {\cal X}_L(M)$, then    $\cal D$ is  called {\bf ${\bf {\cal X}}$- invariant} if for any   $X\in {\cal X}$, the tangent map $T_x\phi^X_t $ sends ${\cal D}_x$ onto $ {\cal D}_{\phi^X_t(x)}$ for all $(x,t)\in \O_X$.  If ${\cal D}$ is $ {\cal X}_{\cal D}-$ invariant we simply say that ${\cal D}$ is {\bf invariant}. \\

{\it Now, we introduce   two  properties of "local triviality" which will play an essential role in the whole paper}:\\

$\bullet$   $\cal D$ is {\bf (locally) lower trivial} (lower trivial for short) if, for each $x\in M$, there exists an open neighbourhood $V$ of $x$,  a smooth map $\Theta:\tilde{\cal D}_x\times  V \ap TM$  (called {\bf  lower trivialization}) such that :
\begin{enumerate}
\item[(i)]  $ \Theta(\tilde{\cal D}_x\times\{y\})\subset {\cal D}_y$ for each $y\in V$

\item[(ii)] for each $y\in V$,  $\Theta_y\equiv \Theta(\;,y):\tilde{\cal D}_x\ap T_yM$ is a continuous operator  and $\Theta_x:\tilde{\cal D}_x\ap T_xM$  is the natural inclusion $i_x$

\item [(iii)] there exists a  continuous operator $\tilde{\Theta}_y: \tilde{\cal D}_x\ap \tilde{\cal D}_y$ such that $i_y\circ \tilde{\Theta}_y=\Theta_y$, $\tilde{\Theta}_y$ is an isomorphism from $\tilde{\cal D}_x$ onto ${\Theta}_y(\tilde{\cal D}_x)$
and  $\tilde{\Theta}_x$ is the identity of $\tilde{\cal D}_x$\\
\end{enumerate}

$\bullet$  ${\cal D}$ is called  {\bf (locally) upper  trivial} (upper trivial for short)  if,  for each $x\in M$, there exists an open neighbourhood $V$ of $x$, a Banach space ${F}$ and  a smooth map $\Psi:F\times  V \ap TM$  (called {\bf upper trivialization})  such that :
\begin{enumerate}
\item[(i)] for each $y\in V$,  $\Psi_y\equiv \Psi(\;,y):F\ap T_yM$ is a continuous operator with $\Psi_y(F)={\cal D}_y$;
\item [(ii)] $\ker \Psi_x$  complemented in $F$;
\item[(iii)] if $F= \ker \Psi_x\oplus S$,  the restriction $\theta_y$  of $\Psi_y$ to $S$  is injective for any $y\in V$;
 \item[(iv)]  ${\Theta}(u,y)=({\theta}_y\circ [{\theta}_x]^{-1}(u), y)$ is  a lower trivialization of ${\cal D}$. 
 \end{enumerate}
In this case the map  $\Theta$ is called the {\bf associated lower trivialization}.

\begin{exes}${}$
\begin{enumerate}
\item  The distribution ${\cal D}_x=x+T(F)$ (where $T: F\ap E$   is a bounded operator as in Example \ref{ex} $n^0$ 2) on $E$ is  lower trivial. This distribution is upper trivial if and only if  $\ker T$ is complemented in $F$
\item Let be $\S$  a  closed topological subset  of a Banach space $E$ and again $T: F\ap E$ a bounded operator with $T(F)\not=E$. We consider the  distributions $\cal D$  and ${\cal D}'$ 
on $E$ defined in the following way:\\
${\cal D}_x=x+T(F)$ if $x\in \S$ and ${\cal D}_x=x+E$ otherwise;\\
${\cal D}'_x=x+E$ if $x\in \S$ and ${\cal D}'_x=x+T(F)$ otherwise.\\
 It is easy to see that $\cal D$ and ${\cal D}'$ are 
 weak distributions on $E$.  Then ${\cal D}$ is  lower trivial but not  upper trivial  and ${\cal D}'$ 
 is neither   lower trivial, nor  upper trivial.
\end{enumerate}
 \end{exes}

\bigskip 

\noindent For an illustration of the property of  lower local triviality  and upper  local triviality in a more large context, we give the following result, which is a generalization of Example  \ref{ex}   
$n^0 3$
\begin{prop}\label{upperloc}${}$\\
Let ${\cal D}: x\ap{\cal D}_x\subset T_xM$ be a field   on $M$ of  normed subspaces. Suppose that for each $x\in M$, there exists an open neighbourhood $V$ of $x$, a Banach space ${F}$ and  a smooth map $\Psi:F\times  V \ap TM$   such that :
%\begin{enumerate}
%\item[(i)]  $ \Psi(F\times\{y\})= {\cal D}_y$ for each $y\in V$
%\item[(ii)] 
for each $y\in V$,  $\Psi_y\equiv \Psi(\;,y):F\ap T_yM$ is a continuous operator such that $\Psi_y(F)={\cal D}_y$
%\end{enumerate}
Then, we have the following properties:
\begin{enumerate}
\item ${\cal D}_x$ has a natural structure of Banach space (again denoted by $\tilde{\cal D}_x$) such that the canonical continuous operator $\hat{\Psi}_x:F/\ker \Psi_x\ap \tilde{\cal D}_x$ is an isometry, in particular, ${\cal D}$ is a weak distribution.
\item    There exists a neighbourhood $W$ of $x$ and, for any $y\in W$, a  continuous surjective  operator $\tilde{\Psi}_y: F\ap \tilde{\cal D}_y$ such that $i_y\circ \tilde{\Psi}_y=\Psi_y$, where $i_y:\tilde{\cal D}_y\ap T_yM$ is the natural inclusion. 
\item   Assume  that $\ker \Psi_x$ is complemented (i.e.  $F=\ker \Psi_x\oplus S$). Then there exists an open neighbourhood $W$ of $x$ such that the restriction $\theta_y$  of $\Psi_y$ to $S$  is injective for any $y\in W$,  and then ${\Theta}(u,y)=({\theta}_y\circ [{\theta}_x]^{-1}(u), y)$ is  {a lower trivialization of ${\cal D}$}. So, ${\cal D}$ is   upper trivial, and {$\Psi$ is an upper  trivialization of $\cal D$}.
\end{enumerate}
\end{prop}
\bigskip

The context of Proposition \ref{upperloc} can be found in the framework of Banach Poisson manifold $(M,\{\;,\;\})$ where  $\Psi :T^*M\ap TM\subset  T^{**}M$ is the canonical morphism associated to the Poisson structure (see for instance  \cite{OdRa1} and \cite{OdRa2}).

\begin{cor}\label{fibre}${}$\\
Let $\pi:{\cal F}\ap M$ be a Banach fiber bundle over $M$ with  typical fiber $F$ and   $\Psi:{\cal F}\ap TM$ a morphism of bundle.  If the kernel of $\Psi$ is complemented in each fiber,  then ${\cal D}$ is an upper   trivial  weak distribution.
\end{cor}

\begin{obs}\label{smoothG}${}$
\begin{enumerate}
\item Recall  the definition of "differentiability"  (resp. "smoothness") of a closed and complemented distribution $\cal D$, introduced  in \cite{ChSt} (resp. \cite{Nu}): there exists some neighbourhood $V$ of $x\in M$ on which $TM$ is trivializable and there exists on $V$ a smooth field $y\ap\hat{\Theta}_y$ of  isomorphisms from $T_xM$ to $T_yM$ so that $\hat{\Theta}_x=Id_{T_xM}$ and $\hat{\Theta}_y[{\cal D}_x]\subset {\cal D}_y$. 

For a closed distribution, we have $\tilde{\cal D}\equiv {\cal D}$ and so $\Theta(u,y)=\hat{\Theta}_y(u)$ $(u,y)\in {\cal D}_x\times V$ is a lower trivialization and so, the property of lower triviality is nothing but a generalization of "differentiability".  Note that, the property of "smoothness" for a closed distributions, introduced in \cite{Gl}  and \cite{Nu}, is  also a  generalization of the property of "differentiability" given in  \cite{ChSt} in the same way.
\item Consider an anchored  bundle on $M$, that is  a Banach bundle  $\pi:{\cal F}\ap M$  over $M$ so that we have    a  bundle morphism $\Psi:{\cal F}\ap TM$ .  In this situation, we are going to see that the distribution ${\cal D}=\Psi({\cal F})$ is "continuous local  lower trivial"  in some sense, which is quite different of our previous definition. At first,
 recall that, from \cite{Mi}, if $p:F\ap F/K$ is the canonical projection of a Banach space $F$ onto the Banach quotient $F/K$ ,  there always exists a continuous selection $\s: F/K\ap F$ such that $\s(\l u)=\l \s(u)$ for any $\l\in \R$ and $u\in  F/K$. However, $\s$ is not linear. Of course, there exists  such a linear selection  $\s $,  if and only if $K$ is complemented in $F$. Come back to  the  situation of an anchored bundle  $\Psi:{\cal F}\ap TM$. We fix  some  $x\in M$,  and choose a continuous selection $\theta_x: {\cal F}_x/ \ker \Psi_x$. We can identify $\tilde{\cal D}_x=\Psi({\cal F}_x)$  with $ {\cal F}_x/ \ker \Psi_x$. Thus, any $u\in {\cal D}_x$ can be written $\Psi\circ \theta_x(u)$.  Choose  any neighbourhood  $U$ of $x$ such that ${\cal F}_{| U}$ is equivalent to ${\cal F}_x\times U$. Via the previous equivalence we consider the map $\theta: {\cal D}_x\times U\ap TM$ defined by:
 $$\theta(u,y)=\Psi(\theta_x(u),y)$$
 Then $\theta$ is continuous and  in general not linear  in the first variable, and  is smooth in the second variable. So, we can consider $\theta$ as a kind of  local lower trivialization which, in general is only continuous  and not linear on the fiber ${\cal D}_x$. 
 
 On the other hand, this construction can be considered as   a kind  of "smoothness" of ${\cal D}$:  via such a (local) map $\theta$,  each $u\in{\cal D}_x$ can be extended to a smooth vector field which is tangent to ${\cal D}$ (compare with \cite{OdRa1} p 35). 
 
 Of course if $\theta _x$  is linear, we get a lower trivialization associated to $\Psi$ (see the proof of Proposition \ref{upperloc}). 
 
 However, in this context,  our criteria of integrability works only when we have a lower trivialization associated to an upper trivialization. For   this reason,  we  have introduced  the property of upper  triviality.
 \item In finite dimension,   the definition of "differentiability"  given in \cite{St} is the same definition  as  that  given in  \cite{ChSt}.  On the other hand, in \cite{Su}, a distribution is called smooth, if there exists a subset $\cal X$ which generates $\cal D$.  In this  context, for any $x\in M$ there exists a neighbourhood $V$ of $x$ and vector fields $X_1,\cdots,X_p$ which are defined and linearly independent on $V$ such that $\{X_1(x),\cdots,X_p(x)\}$ is a basis of ${\cal D}_x$ so $\cal D$ is lower trivial.
 On the other hand, for  any $x\in M$,   denote by ${\cal X}_{{\cal D}_x}$ the module of germs vector fields at $x$ which are tangent to $\cal D$. If   ${\cal X}_{{\cal D}_x}$ is finitely generated, then $\cal D$ is  upper  trivial.
  
\end{enumerate}
\end{obs}
\bigskip

We end this section with the proof of  Proposition \ref{upperloc} and its Corollary. For this proof, we need the following Lemma which will be also used later:

\begin{lem}\label{inject}${}$
\begin{enumerate}
\item Consider  two Banach spaces ${E_1}$ and  $E_2$ and $i:E_1\ap E_2$ an injective continuous  operator.  Let $y\ap \Theta_y$ be a  smooth field of continuous operators  of $L(E_1,E_2)$ on an open neighbourhood $V$ of $x\in E_1$ such that $\Theta_x=i$. Then there exists  a neighbourhood  $W$ of $x$ in $V$ such that for each $y\in W$, $\Theta_y$ is an injective operator.
\item Let $f:U \ap V$ be a $C^1$ map from two open sets $U$ and $V$ in Banach spaces $E_1$ and $E_2$ respectively, such that $T_uf$ is injective at $u\in U$.  Then there exists an open neighbourhood $W$ of $u$ in $U$ such that the restriction of $f$ to $W$ is injective.
\end{enumerate}
\end{lem}

\noindent\begin{proof}\so{\it Proof of Lemma \ref{inject}}${}$\\
There exists an open ball $B(x,r)$ included in $V$ such that $||\Theta_{y}-\Theta_x||\leq K||y-x]]$ for any $y\in B(x,r)$. We can suppose that $r<1$. Assume that the conclusion of  Lemma \ref{inject} (1) is not true. So, for each $n\in \N^*$, there exists $x_n\in B(x,r/n)$ and $h_n\in E_1$ with $||h_n||=1$ such that $\Theta_{x_n}(h_n)=0$. We have of course:
\begin{eqnarray}\label{nulsuite}
<\a,\Theta_{x_n}(h_n)>=0
\end{eqnarray}
for all $\a\in E_2^*$.\\
It follows that we have:
\begin{eqnarray}\label{inegx}
|<\a ,\Theta_x(h_n)>|=|<\a,(\Theta_x-\Theta_{x_n})(h_n)>|\leq\dis\frac{ K||\a||}{n}
\end{eqnarray}
On the other hand, $\Theta_x=i$ is a continuous bijective operator from the Banach space $E_1$ onto the normed subspace $F=i(E_1)$ in $E_2$. So, the transpose operator $i^*\in L(F^*,E_1^*)$ is a monomorphism\footnote{ an operator $T$  between two Banach space $E$ and $F$ is a monomorphism if we have\\ $\inf \{|| T(u)||_F;\; ||u||_E=1\}\geq k>0$}
with a dense range (see \cite{Ha}, \cite{HaMb}). According to Hahn-Banach Theorem, there exists $\b_n\in E_1^*$ such that $<\b_n,h_n>=1$ with $||\b_n||=1$. From the density of $i^*(F^*)$, there exists $\a_n\in F^*$ such that $||\b_n-i^*(\a_n)||<\dis\frac{1}{4}$, i.e. such that $\dis\frac{3}{4}\leq ||i^*(\a_n)||\leq \dis\frac{5}{4}$ . From these inequalities we get:\\

$\bullet$ $|<i^*(\a_n),h_n>-1|=|<i^*(\a_n)-\b_n,h_n>| \leq \dis\frac{1}{4}$ \\
so we have $|<i^*(\a_n),h_n>|\geq \dis\frac{3}{4}$ \\

$\bullet$ as  $||i^*(\a_n)||\leq \dis\frac{5}{4}$ and as $i^*$ is a monomorphism, we have $||i^*(\a_n)||\geq k||\a_n||$ for some $k>0$ and finally we get $||\a_n||\leq \dis\frac{5}{4k}$. \\

On the other hand we can write:
\begin{eqnarray}\label{min}
|<i^*(\a_n),h_n>|=|<\a_n,i(h_n)>|\geq \dis\frac{3}{4}
\end{eqnarray}
for any $n$. \\
From Hahn-Banach Theorem, we obtain the same relation (\ref{min})  with $\a_n\in E_2^*$.  But from (\ref{inegx}) we get:

$|<\a_n ,\Theta_x(h_n)>|\leq \dis\frac{ K||\a||}{n}$ and so  $|<\a_n,i(h_n)>| \leq \dis\frac{K}{n}\frac{5}{4k}$ \\for any $n$.

\noindent So we get a contradiction with (\ref{min}) for $n$  large enough. So we have completed the proof of the part (1).\\

 Let  $f:U \ap V$ be a map of class $C^1$.  As the problem is local,  without loss of generality, we can suppose that $U$ is an open ball  of center $0\in E_1$. As $f$ is of class  $C^1$, there exists an open ball $B(0,r)$ such that
\begin{eqnarray}\label{lipDphi}
||T_uf -T_vf ||\leq K ||u-v|| \textrm{ for } u,v \in B(0,r)
\end{eqnarray}

\noindent Moreover,  we can choose $r$ so  that $r<1$.\\

Suppose that $f$ is not locally injective around $0$. Given any pair $(u,v)\in[ B(0,r)]^2$ such that $u\not =v$ but $f(u)=f(v)$,
we set $h=v-u$.  For any $\a \in E_2^*$ we consider the smooth curve $c_\a:[0,1]\ap \R$ defined by:
$$c_\a(t)=<\a ,f(u+th)-f(u)>$$
Of course we have  $\dot{c}_{\a}(t)=<\a, T_{u+th}f(h)>$.\\

\noindent Denote by $]u,v[$ the set of points $\{w=u+th, t\in ]0,1[\}$. As we have $c_\a(0)=c_\a(1)=0$, from Rolle's Theorem, there exists $u_\a\in ]u,v[$ such that
\begin{eqnarray}\label{AF}
<\a,T_{u_\alpha}f(h)>=0
\end{eqnarray}
Replacing $h$ by $\ds\frac{h}{||h||}$, we can suppose in (\ref{AF}) that $||h||=1$.

From our assumption it follows that,  for each $n\in \N^*$, there exists $u_n$ and $v_n$ in $B(x,r/n)$ so that $u_n\not =v_n$ but with $f(u_n)=f(v_n)$. So from the previous argument, for any $\a \in E_2^*$, we have
\begin{eqnarray}\label{AFn}
<\a,T_{u_{\alpha,n}}f(h_n)>=0
\end{eqnarray}
for some $u_{\a, n}\in]u_n,v_n[$ and with $h_n=\ds\frac{v_n-u_n}{||v_n-u_n||}$

From (\ref{lipDphi}) and (\ref{AFn}), we get
\begin{eqnarray}\label{inega}
|<\a,T_0f(h_n)>|=|<\a,[T_0f-T_{u_{\a,n}}f](h_n)>| \leq  ||\a||.\dis\frac{Kr}{n}<||\a||\dis\frac{K}{n}.
\end{eqnarray}
for any $\a\in E_2^*$.

Now, we can  use the same argument as in part (1) and we get again  a contradiction.

\end{proof}

\noindent \begin{proof}\so{\it Proof of Proposition \ref{upperloc}} ${}$\\

At first, for any $x\in M$,  we have a natural Banach structure on ${\cal D}_x$ (again denoted by $\tilde{\cal D}_x$) such that the natural morphism  $\tilde{\Psi}_x  :F/\ker \Psi_x\ap \tilde{\cal D}_x$ is an isometry. On the other hand, take a local trivialization of $TM$ on a neighbourhood $W$ of $x$; so we have $TM\equiv E\times W$. In this context,   on $W$,  $\Psi$ can be identified with a smooth field of continuous operators
$\Psi_y:F\ap E$ such that  ${\cal D}_y=\Psi_y(F)\times\{y\}\subset E\times \{y\}\equiv T_yM$. Let us consider the following commutative diagram:

$$\xymatrix {
    F  \ar[r]^{^q\;\;\; \;\;\;\;}\ar[d]_{_{\Psi_y}} & {F/\ker{ \Psi_y}} \ar[ld]^{\hat{\Psi}_y} \\
   {\cal D}_y
  }$$
where $q$ is the natural projection and ${\hat{\Psi}_y} $ is the natural bijection induced by ${\Psi_y}$. So, if we consider the  Banach structure $\tilde{\cal D}_y$, we get a continuous operator $\tilde{\Psi}_y={\hat{\Psi}_y}\circ q:F\ap \tilde{\cal D}_y$ so that $\Psi_y=i_y\circ \tilde{\Psi}$. \\

Assume that $F=\ker \Psi_x\oplus S$, for some Banach space $S\subset F$. Let ${\theta}_y$ be the restriction to $S$ of $\Psi_y$ for any $y\in W$. Clearly, ${\theta} (u,y)=(\theta_y(u),y)$ defines a smooth map from $S\times W$ into $ E\times V\equiv TM$ and $\theta_y:S\times\{x\}\ap E\times \{x\}\equiv T_xM$ is a continuous operator whose image is contained in ${\cal D}_y$. \\

On the other hand, let  $\tilde{\theta}_y$ be the restriction of  $\tilde{\Psi}_y$ to $S$, then, $\tilde{\theta}_y$ is a continuous operator from $S$ to $\tilde{\cal D}_y$ so that $\theta_y=i_y\circ \tilde{\theta}_y$ for any $y\in W$. Of course, $\tilde{\theta}_x :S\ap \tilde{\cal D}_x$ is an isometry and, in particular, it is an isomorphism.  As, $\theta_x$ is injective, according to Lemma  \ref{inject}, without loss of generality, we can suppose that $\theta_y$ is injective for any $y\in W$. It follows that $\tilde{\theta}_y$  is a continuous injective operator from $S$  into $\tilde{\cal D}_y$. As $\theta_y$ is injective,
we have $\ker\Psi_y\cap S=\{0\}$. It follows that $q_1=q_{|S}$ is an isomorphism onto $q(S)\subset F/\ker\Psi_y$. Of course the restriction $q_2$ of the isomorphism $\hat{\Psi}_y: F/\ker \Psi_y\ap \tilde{\cal D}_y$ to $q(S)$ is an isomorphism  onto $\tilde{\theta}_y(S)$  such that
$\tilde{\theta}_y=q_2\circ q_1$. So $\tilde{\theta}_y$ is an isomorphism of $S$ onto $\tilde{\theta}_y(S)$.\\

Finally, the map
$${\Theta}:\tilde{\cal D}_x\times W\ap E\times W\equiv TM$$
defined by   ${\Theta}(u,y)=({\theta}_y\circ [{\theta}_x]^{-1}(u), y)$ is clearly a lower trivialization of ${\cal D}$

\end{proof}

\noindent\begin{proof}\so{Proof of Corollary \ref{fibre}}${}$\\

Given $x\in M$ there exists a local trivialization of ${\cal F}$ on an open set $V$ around $x$. So we can identify $\cal F$ with $F\times V$ on $V$. In this context, in restriction to $V$, the morphism $\Psi$ can be identified,  as a map $\Psi: F\times V\ap TM$ which satisfies assumption (i) and (ii) of Proposition \ref{upperloc}

\end{proof}

%%%%%%%%%%%%%%%%%%
\subsection{Results}\label{resultats}
%%%%%%%%%%%%%%%%%
%\begin{rem}\label{lowsection}${}$
 Let ${\cal D}$   be a lower  trivial   distribution on $M$. For any $x\in M$ and any  lower trivialization  \\$\Theta :\tilde{\cal D}_x\times  V\ap TM$
and any $u\in \tilde{\cal D}_x$  ,we  consider
  \begin{eqnarray}\label{localsec}
X(z)=\Theta(u,z)
\end{eqnarray}
Of course we also have $X(z)=i_z\circ\tilde{\Theta}(u,z)$ where $\tilde{\Theta}(u,z)=\tilde{\Theta}_z(u)$ and  $X$ is a local vector field on $M$ tangent to ${\cal D}$  whose domain is $V$. Moreover, the set of all such local vector fields spans ${\cal D}$.
%\end{rem}

A {\bf  lower (local) section} of a lower   trivial weak  distribution ${\cal D}$ is  a map of type (\ref{localsec})  for  any lower trivialization $\Theta$ any  $u\in \tilde{\cal D}_x$  and any $x\in M$. Note that  the domain of  all lower  section defined by  such a lower trivialization $\Theta:\tilde{\cal D}_x\times V\ap TM$ is the open set $V$. The {\bf set of such lower sections} will be denoted by $\bf {\cal X}^-_{\cal D}$. \\

The following Proposition gives a relation between integral manifolds and  ${\cal X}^-_{\cal D}-$invariant  weak distributions:

\begin{prop}\label{lem1}${}$\\
If  a lower  trivial weak  distribution ${\cal D}$  (resp.  lower  trivial  closed distribution)  is integrable, then ${\cal D}$  is $ {\cal X}^-_{\cal D} -$invariant (resp.  $ {\cal X}_{\cal D} -$invariant)
  \end{prop}

In this context, we obtain the following  version of Stefan-Sussmann Theorem:

\begin{theor}\label{thI}
Let ${\cal D}$ be a lower trivial weak  distribution on a Banach manifold $M$.
\begin{enumerate}
\item ${\cal D}$ is integrable if and only if it is ${\cal X}^-_{\cal D}$-invariant.
\item if $\cal D$ is integrable, on $M$,
consider the binary relation
$$x{\cal R}y \textrm { iff there exists an integral manifold }( N,f) \textrm{  of  } {\cal D} \textrm{ such that } x ,y \in  f(N).$$
Then $\cal R$ is an equivalence relation and the equivalence class $L(x)$ of $x$ has a natural structure of connected Banach manifold modelled on $\tilde{\cal D}_x$.\\
Moreover  $(L(x), i_{L(x)})$, is a maximal integral manifold of ${\cal D}$  in the following sense: for any integral manifold $(N,f)$ of $\cal D$, such that $f(N)\cap L(x)$ is not empty  then $f(N)\subset L(x)$.
%\end{enumerate}
\end{enumerate}
\end{theor}

Taking into account Remark \ref{closed}, the property of  lower triviality of a  weak distribution corresponds to the usual assumptions on the distribution that we find in \cite{ChSt}, \cite{Nu}, \cite{Gl}. When ${\cal D}_x$ is closed (resp. complemented) in $T_xM$ the following Corollary of Theorem 1 gives exactly the main result of integrability of distributions  we can find in  \cite{Gl} (resp. \cite{ChSt}, \cite{Nu}):

\begin{cor}${}$ \\ For a lower trivial closed distribution the following propositions are equivalent:
\begin{enumerate}
\item[(i)] $\cal D$ is integrable;
\item[(ii)] $\cal D$ is ${\cal X}_{\cal D}-$invariant;
\item[(iii)] $\cal D$ is ${\cal X}^-_{\cal D}-$invariant.
\end{enumerate}
\end{cor}

\bigskip
We end this section with the \so{\it proof of Proposition \ref{lem1}} :

\noindent \begin{proof}

\noindent Consider  a lower section $X(y)=\Theta(u,y)$, associated to a lower trivialization    $\Theta:\tilde{\cal D}_x\times  V\ap TM$. So Dom$(X)=V$. Fix such such a lower section $X$ and the associated lower trivialization $\Theta$. Denote by $\D$ (resp. $\tilde{\D}$) the subspace ${\cal D}_x$ of $T_xM\equiv E$ (resp. the Banach space $\tilde{\cal D}_x$). Let be $\D_y=\Theta_y(\D)$ and $\tilde{\D}_y$ the natural Banach structure induced by $\tilde{\cal D}_y$.\\

Given any $z\in V$, the map $\Theta'_y=\Theta_y\circ[\tilde{\Theta}_z]^{-1}$ is a smooth field  of continuous operators from $\tilde{\D}_z$ into $T_yM\equiv E\times\{y\}$ and moreover, $\tilde{\Theta'}_y=\tilde{\Theta}_y\circ[\tilde{\Theta}_z]^{-1}$ is an isomorphism between $\tilde{\D}_z$  and $\tilde{\D}_y$. Of course, if $v=\tilde{\Theta}_z(u)$, we have  $X(y)=\Theta'_y(v)$.\\

Let $f:N\ap M$ be an integral manifold of $\cal D$  passing through some  $z\in V$. Then,  $N$ is a Banach manifold  modelled on the Banach space $\tilde{G}=\tilde{\D}_z$. For any open neighbourhood $U$ of $z$ the set  $\tilde{U}=f^{-1}(U)$ is an open neigbourhood of $\tilde{z}=f^{-1}(z)$. According to Remark \ref{weakN}, without loss of generality,  we may assume that $N$ is an open set in $\tilde{G}$ and $M$ is an open set in $E$. In these identifications, $f$ is the natural inclusion $i_N$ of $N$ in $M$, that is the restriction to $N$ of the natural inclusion $i:\tilde{G}\ap E$. In this context, on $i(N)\subset M$, $y\ap \Theta'_y$ is  a smooth  field of continuous linear operator from $\tilde{\D}_z\subset \tilde{G}$  into ${\cal D}_y\equiv G\times\{y\}$. Moreover,  $\tilde{\Theta'}_y$ is an isomorphism  between $\tilde{\D}_z\subset \tilde{G}$ and $\tilde{\D}_y\subset \tilde{G}\times\{y\}$ for any  $y\in i(N)$.\\

\begin{lem}\label{Xdiff}${}$\\
With the previous notations, the map $y\mapsto \tilde{\Theta'}_y$ from $N$ to $L(\tilde{\D}_z,\tilde{G})$ is smooth (for the topology induced by $\tilde{\D}_z$ on $N$).
\end{lem}

From Lemma \ref{Xdiff}, $\tilde{Y}=\tilde{\Theta'}_y(v)$ is a smooth vector field on the  Banach manifold $N$, and, moreover,
  and we have  $X(i(y))=T_yi[\tilde{Y})(y)]=(i_*Y)(y)$ on $i(N)$. So  the flow $\phi^X_t$ satisfies the relation
$$\phi^X_t\circ i=i\circ \phi^{\tilde{Y}}_t$$
 on a small neighbourhood  $W$ of $ z$ and for all $t$ such that $\phi^{\tilde{Y}}_t$ is defined on $N$. Of course for any $y\in W$  and  $t$  such that $\phi^{\tilde{Y}}_t(y)$ is defined, we have
 $$T_y\phi^X_t({\cal D}_y)=T_y\phi^X_t(i[T_yN])=i\circ T_y\phi^{\tilde{Y}}_t(T_yN)=i[T_{\phi^{\tilde{Y}}_t(y)}N]=i[\tilde{\cal D}_{\phi^{\tilde{Y}}_t(y)}]={\cal D}_{\phi^X_t(y)}.$$

 Now, consider  any $(z,t)\in \O_X$. Denote by $]\a_z,\b_z[$ the maximal interval on which $\Phi^X_t(z)$ is defined. Given any $\t\in]\a_z,\b_z[$, consider the integral curve $\g(t)=\phi^X_t(z)$ for $t\in [0,\t]$.  By compactness  of $[0,\t]$ there exists a finite number of integral manifolds $(N_1,f_1),\cdots, (N_r,f_r)$ so that $\g([0,\t])$ is contains in $\dis\cup_{i=1}^r f_i(N_i)$. Using the previous argument, by induction, we obtain:

 $$T_z\phi^X_\t({\cal D}_z)={\cal D}_{\phi^X_\t(z)}.$$
 We deduce that integrability implies $ {\cal X}^-_{\cal D} -$invariance.\\

Now, if moreover ${\cal D}$ is closed, given an integral manifold $f:N\ap M$ and any local section $X$ of ${\cal D}$  whose domain intersects $f(N)$, then  $X$ induces, by restriction on $f(N)$,   a smooth vector fields on $N$.  So the same arguments used last part of in the previous  proof   works  too. ( see \cite{Gl}).\\
 \end{proof}

\noindent \begin{proof}\so{\it Proof of Lemma \ref{Xdiff}}${}$\\

 From convenient analysis (see \cite{KrMi}), recall that for a map $f$ from an open set $U$ in a Banach space $E_1$ to a Banach space $E_2$ we have the equivalent following properties
 \begin{enumerate}
 \item[(i)] $f$ is smooth;
 \item [(ii)] for any smooth curve $c:\R\ap U$ the map  $t\mapsto f\circ c(t)$ is smooth;
 \item [(iii)] the map $t\mapsto <\a, f\circ c(t)>$ is smooth for any $\a\in E^*_2$;\\
 \item [(iv)] for any smooth curve $c:\R^2\ap U$, all partial derivatives of $f\circ c\;$ exist and are locally bounded. \\
  \end{enumerate}
  \bigskip
 Fix some $v\in \tilde{\D}_z $. Note that, for any $\a\in \tilde{G}^*$  we have
\begin{eqnarray}\label{aX}
<\a,\tilde{\Theta'}_y(v)>=<[\tilde{\Theta'}_y]^*(\a),v>
\end{eqnarray}

If $i:\tilde{G}\ap G$ is the natural inclusion, we have $[\Theta'_y]^*=[\tilde{\Theta'}_y]^*\circ i^*$.\\

For $y\in i(N)\subset \tilde{\D}_z$  and $\a\in G^*$ fixed, we  consider the map

$ h(y)=[\Theta'_y]^*(\a)=[\tilde{\Theta'}_y]^*(i^*\a))$ \\
Clearly, $h$ is a smooth map from  the open $i(N)$ in the normed space $\D_z\subset E$ to  the Banach $[\tilde{\D}_z]^*$. Take any smooth curve $c:\R\ap N\subset\tilde{\D}_z$. As the inclusion of $\tilde{\D}_z$ into $\D_z$ is linear continuous, $c$ is also a smooth map from $\R$ to $N\subset \D_z$, the map $h\circ c$ is a smooth map from $\R$ to $[\tilde{\D}_z]^*$. We conclude that $h$ is a smooth map from $N\subset \tilde{\D}_z$ to $[\tilde{\D}_z]^*$

So from (\ref{aX}), we see that  the map $y\ap <i^*\a,\tilde{\Theta'}_y(v)>$ is a smooth from $N\subset \tilde{\D}_z$ to $\R$, for any $\a\in G^*$.   As $i^*({G}^*)$ is dense in $\tilde{G}^*$, given any $\b\in \tilde{G}^*$ there exists a sequence $\a_n\in {G}^*$ so that $i^*(\a_n)$ converges to $\b$ in $\tilde{G}^*$.  For simplicity, we set
$g(y)= \tilde{\Theta'}_y(v)$. Consider any smooth curve $c:\R\ap N\subset \tilde{\D}_z$.

Now on any compact $K\subset \R$, and for any $p\in \N$ we have:

$|<\b,(g\circ c)^{(p)}(t)>-<i^*(\a_n),(g\circ c)^{(p)}(t)>|\leq ||\b -i^*(\a_n)||\dis \sup_{t\in K}|(g\circ c)^{(p)}(t)|$

So the map  $<i^*\a_n,(g\circ c)^{(p)}>$ converges uniformly to $<\b,(g\circ c)^{(p)}>$ on $K$. It follows that $<\b,g\circ c>$  is a smooth map for any $\b\in\tilde{G}^*$.  On one hand, we have proved that  the map $y\ap \tilde{\Theta'}_y(v)$is smooth for any $v\in \tilde{\D}_z$. On the other hand, we know that  $ \tilde{\Theta'}_y$ is a continuous operator from $\tilde{\D}_z$  to $\tilde{G}$. I It follows from (iv)  that the map $y\ap  \tilde{\Theta'}_y$ is a smooth map from $N\subset \tilde{\D}_z$ into $L(\tilde{\D}_z,\tilde{G})$.\\
\end{proof}

%%%%%%%%%%%%%%%%%%%%%%%%%%%%%%%%%%%%%%%%%%%%%%%%%%%
\subsection{Proof of Theorem 1}\label{thI}
%%%%%%%%%%%%%%%%%%%%%%%%%%%%%%%%%%%%%%%%%%%%%%%%%%
{\bf Proof of Part (1)}${}$\\
%%%%%%%%%%%%%%%%%%%%%%%%%%%%%%%%%%%%%%%%%%%%%%%%%%
At first, according to the Proposition \ref{lem1}, integrability implies $ {\cal X}^-_{\cal D} -$invariance.\\
So we have to prove the converse.
In fact, this proof is an adaption of arguments of Chillingworth and Stefan used in \cite{ChSt}.\\

Given $x\in M$, we may assume that  $M$ is an open set of $E$ and $TM\equiv E\times M$. We denote by $\D$ (resp. $\tilde{\D}$ )
  the normed space (resp. the  Banach space)  ${\cal D}_x$ (resp. $\tilde{\cal D}_x$).  From the property of lower  triviality,  after
  restricting this open if necessary, we have  a smooth  fields  $y\ap \Theta_y$  of  continuous operators from $\tilde{\D}$ to $E$. Consider the family $\{X_u(y)= \Theta_y(u),\; u\in \tilde{\D}\}$ of smooth vector fields  on $M$. By standard argument (see \cite{ChSt} proof of Corollary 4.2), we can choose an open ball $B(0,r)\subset\tilde{ D}$ so that the flow $\phi^{X_u}_t$ is defined on an open neighbourhood $W$ of $x$ for all $|t|\leq 1$.
   We set  $\Phi(t,y,u)= \phi^{X_u}_t(y)$, $t\in [0,1]$, $y\in W$ and $u\in B\equiv B(0,r)\subset \tilde{\D}$

   \begin{lem}\label{invariance}${}$\\ For any smooth map  $\Phi:\R\times W\times B\ap E$  we denote by $D_t\Phi(t,y,u)$  (resp. $D_y\Phi(t,y,u)$, resp. $D_u\Phi(t,y,u)$ the partial derivative of $\Phi$ according to the first (resp.  the second (resp. the third) variable, at point $(t,y,u)\in \R\times,W\times B$. \\
With these notations, $u\ap \Phi(t,y,u)$ is smooth . Moreover   assume that  $T_x\phi^{X_u}_t[{\cal D}_x]={\cal D}_{\phi^{X_u}_t(x)}$ for all $t$ such that $(x,t)\in \O_{X_u}$ and all $u\in B$, then we have:
   \begin{eqnarray}\label{inclusion}
D_u\phi(t,x,u)(\D)\subset {\cal D}_{x(t)}
\end{eqnarray}where $x(t)=\phi(t,y,u)$.
\end{lem}

\noindent \begin{proof}\so{\it Proof}${}$

\noindent At first, we fix $y\in W$ and   $u\in B$, and we  set :\\
$y(t)=\phi(t,y,u)$ (the integral curve of $X_u$ though $y$);\\
$X(t,y,u)=X_u(y(t,u))$;\\
 $A(t)=D_yX(t,y,u)$;\\
  $B(t)=D_uX(t,y,u)$. \\

  \noindent Of course, $A$ (resp. $B$) is a smooth field of operators  in $L(E,E)$ (resp. $L(\tilde{\D},E)$). In fact, we have  $B(t)=\Theta_{y(t)}$.  So, in the Banach space  $L(\tilde{\D},E)$), the linear differential equation
$$\dot{\S}= A\circ \S+B$$
as an unique solution with initial condition $U(0)=0$ given by

\begin{eqnarray}\label{U(tu)}
\S-t(u)=\ds\G_t\int_0^t(\G_s)^{-1}\circ \Theta_{y(s)}ds
\end{eqnarray}

where $\G_s$ is the solution of the differential equation
$$\dot{\G}=A\circ \G$$
with initial condition $\G_0=Id_E$

From (10.7.3) and (10.7.4) of \cite {Di}, we obtain that $\phi$ is smooth in the third variable and we have

\begin{eqnarray}\label{diffphi}
D_u\phi(t,y,u)=\S-t(u).
\end{eqnarray}

We now look for  the integral curve $x(t)$ through $x$.  In this case,  $\G_s$ is  in fact the  $t\ap \phi^{X_u}_t(x)$  (see\cite{Di} (10.8.5)), taking in account  our assumption of invariance by $\phi^{X_u}_t(x)$, we have:
\begin{eqnarray}\label{inv}
\G_s({\cal D}_x)={\cal D}_{x(s)}
\end{eqnarray}
On the other hand, from the assumption of lower triviality,  we have $\Theta_{x(s)}({\cal D}_x)\subset {\cal D}_{x(s)}$. So,  we get

 $$(\G_s)^{-1}\circ \Theta_{x(s)}({\cal D}_x)\subset  {\cal D}_{x}$$
 and moreover by integration we also have
  $$\ds\int_0^t(\G_s)^{-1}\circ \Theta_{x(s)}({\cal D}_x)\subset  {\cal D}_{x}$$
Finally, using   (\ref{inv}) and (\ref{U(tu)}), we obtain the announced result.\\

   \end{proof}

\noindent We are now in situation to give a sufficient condition of the existence of an integral manifold through $x\in M$:

\begin{prop}\label{immersion}${}$\\
Consider the map:
\begin{eqnarray}\label{diffphi}
\Phi:B\ap M \textrm{ defined by } \Phi(u)=\Phi(1,x,u)= \phi^{X_u}_1(x), \textrm{ for } u\in B\equiv B(0,r)\subset \tilde{\D}.
\end{eqnarray}

 There exists $\d>0$ such that $\Phi:B(0,\d) \ap M$ is a weak submanifold of $M$. Moreover,  if  we have $T_x\phi^{X_u}_t[{\cal D}_x]={\cal D}_{\phi^{X_u}_t(x)}$ for all $t$ such that $(x,t)\in \O_{X_u}$ and all $u\in B$, then, for $\d>0$ small enough,  $(B(0,\d),\Phi)$ is an integral manifold  of $\cal D$ through $x$
 \end{prop}

\noindent It is clear that  Proposition \ref{immersion}, {\bf ends the proof of part (1) of  Theorem 1}.\\

 We now end this subsection with the proof of  the previous  Proposition.\\

\noindent  \begin{proof} \so{\it Proof}${}$\\

 According to Lemma \ref{inject}, it follows that, for $\d>0$ small enough, $(B(0,\d),\Phi)$ is a weak submanifold of $M$.\\

Assume now that  $T_x\phi^{X_u}_t[{\cal D}_x]={\cal D}_{\phi^{X_u}_t(x)}$ for all $t$ such that $(x,t)\in \O_{X_u}$ and all $u\in B$. From Lemma \ref{invariance}, for any $u\in B\subset \tilde{\D}$, we have:
$$D_u{\Phi}(\tilde{\D})\subset {\cal D}_{{\Phi}(u)}.$$
So, it follows that
$$T_u{\Phi}(\tilde{\D})\subset {\cal D}_{{\Phi}(u)}$$ for all $u\in {B}$.\\
Now, according to  the assumption of invariance, we have
\begin{eqnarray}\label{surj}
[T_x\phi^{X_u}_1]^{-1}\circ T_u {\Phi}(\tilde{\D})\subset[T_x\phi^{X_u}_1]^{-1}({\cal D}_{{\Phi}(u)})={\cal D}_x\equiv {F}
\end{eqnarray}

We set $\L_u=[T_x\phi^{X_u}_1]^{-1}\circ T_u{\Phi}$ for $u\in {B}$. In particular, $\L_u$ is a continuous operator from the Banach  space $\tilde{\D}$ to the normed space $\D$. The part (1) will be a consequence of the following Lemma:

\begin{lem}\label{opensurj}${}$\\
Let be $E_1$ (resp. $E_2$) a Banach space (resp. a normed space). Suppose that the set $L_s(E_1,E_2)$ of surjective operators in $L(E_1,E_2)$ is non empty. Then, $L_s(E_1,E_2)$ is an open set.
\end{lem}

\noindent \begin{proof}\so{\it Proof}${}$\\

The first part of this proof is an adaptation of an argument which  can be  found in \cite{QuZu}.\\

Recall that an operator, $T\in L(E_1,E_2)$ is almost open, if for any open ball $B(0,r)$ in $E_2$, there exists an open ball $\tilde{B}(0,\r)\subset E_1$ such that :
$${B(0,r)}\subset \overline{T(\tilde{B}(0,\r))}$$

Given  $\a\in ]0,1[$, there exists $\r>0$ such that, for any $y\in B(0,1)$ we can find $x_1\in \tilde{B}(0,\r)$ such that
$|| y-T(x_1)||\leq \a$. So, $\dis\frac{1}{\a}|| y-T(x_1)||\leq 1$,
and then, there exists $x_2\in \tilde{B}(0,\r)$ such that
\begin{center}
$||\dis\frac{1}{\a}(y-T(x_1))-T(x_2)||\leq \a$ i. e.
$||y-T(x_1)-\a T(x_2)||\leq \a^2$.
\end{center}

By induction, we can build a sequence $(x_n)$ such that $x_n\in  \tilde{B}(0,\r)$ and also
$$||y-T(x_1+\a x_2+\cdots+\a^{n-1}x_n)||\leq \a^n.$$
In the Banach space $E_1$, the series of general term $||\a^{n-1}x_n||$ converges. So, there exists $z\in E_1$ such that $z=\ds\sum_{n=1}^\infty \a^{n-1}x_n$, with $||z||\leq \ds\frac{\r}{1-\a}$ and, of course, $y=T(z)$. It follows that $T$ must be surjective. On the other hand, the set of almost open operator in $L(E_1,E_2)$ is an open set  (see \cite{Ha}, \cite{HaMb}), so the Lemma is proved.

\end{proof}

Coming back to the proof of part(1), the map  $T_0\Phi$ is the inclusion map of $\tilde{F}$ in $F$  and $[T_x\phi^{X_0}_1]=Id_E$  so $\L_0$ is surjective. From Lemma \ref{opensurj}, for $\d>0$ small enough, $\L_u$ is surjective for all $u\in B(0,\d)$; in particular, we get an equality
$$[T_x\phi^{X_u}_1]^{-1}\circ T_u {\Phi}(\tilde{F})=[T_x\phi^{X_u}_1]^{-1}({\cal D}_{{\Phi}(u)})$$
 in (\ref{surj}) which ends the proof of Proposition \ref{immersion}.\\

\end{proof}
%%%%%%%%%%%%%%%%%%%%%%%%%%%%%%%%%%%%%%%%%%%%%%%%%%%
{\bf Proof of Part (2) }${}$\\
%%%%%%%%%%%%%%%%%%%%%%%%%%%%%%%%%%%%%%%%%%%%%%%%%%%%
In this subsection, we will use the  notations introduced in the previous one. In particular, for any $x\in M$, we associate an integral manifold  $(B(0,\d),\Phi)$ build in  Proposition \ref{immersion}. Such an integral manifold will be called a {\bf slice} through $x$. \\
At first, we must prove that the relation $\cal R$ is transitive. This fact is a direct consequence of the following Lemma:
\begin{lem}\label{interplaq}${}$
\begin{enumerate}
\item Given any integral manifold $(N,f)$ of ${\cal D}$ through $x\in M$, there exists a slice $(B(0,\d),\Phi)$ such that $\Phi(0)=x$ and  $f^{-1}[\Phi(B(0,\D)]$ is an open set in $N$
\item For any two integral manifolds $(N,f)$ and $(N',f')$ through $x\in M$, then $f^{-1}[f(N)\cap f'(N')]$ (resp. ${f'}^{-1}[f(N)\cap f'(N')]$) is open in $N$ (resp. $N'$). Moreover, $L=f(N)\cup f'(N')\subset M$ has a natural structure of Banach manifold modelled on $\tilde{\cal D}_x$ and $(L,i_L)$is an integral manifold of $\cal D$ through $x$, where  $i_L$ is the natural  inclusion of $L$ in $M$.
\end{enumerate}
\end{lem}

\noindent \begin{proof}\so{\it Proof}${}$\\

\noindent We fix any $x\in f(N)$. Note that  $N$ is a connected Banach manifold modelled on  $\tilde{\D}\equiv \tilde{\cal D}_x$. As the problem is local,   according to Remark \ref{weakN}, we can assume that $N$ is an open subset of $\tilde{\D}$,  $M$ is an open subset of $E\equiv T_xM$ and $f$ is the natural inclusion $i$ of $\tilde{\D}$ into $E$ (restricted to $N$). Consider a lower trivialization $\Theta:\tilde{\D}\times V\ap M$ around $x$. Given any $u\in \tilde{\D}$, according to the  arguments  used in the proof of Proposition \ref{lem1}, (with $\Theta $ instead of $\Theta'$),  we get that the restriction of  $X_u=\Theta(u,\;)$ to $i(N)$ induces a vector field  $\tilde{Y_u}$  on $i(N)$  relative to  its natural Banach manifold structure.  It follows that the integral curve $t\ap \Phi^{X_u}_t(x)$ of $X_u$ through $x$ lies in $ i(N)$. So, for $\d$ small enough, $\Phi[B(0,\d)]$ is contained in $i(N)\subset i(\tilde{\D})\subset E$. But  as sets, we have $i(\tilde{\D})=\tilde{\D}=\D$. So using  the same arguments used in the proof of part (1) of Theorem 1,  but  in the Banach space $\tilde{\D}$, we can prove that $\Phi$ is a local diffeomorphism of $B(0,\d)$ into $N$ for $\d$ small enough. In particular $L=\Phi[B(0,\d)]$ is an open subset  for the topology of the Banach structure on $i(N)$, which ends the proof of part (1).\\

Let be $(N,f)$ and $(N',f')$ integral manifolds through $x\in M$. Note that $N$ and $N'$ are connected Banach manifold modelled on $\tilde{\D}\equiv \tilde{\cal D}_x$. Applying  part (1) for any $z\in f(N)\cap f'(N')$  to the integral manifold $(N,f)$ (resp. $(N',f')$ we obtain that  $f^{-1}[f(N)\cap f'(N')]$ (resp. ${f'}^{-1}[f(N)\cap f'(N')]$) is open in $N$ (resp. $N'$).\\
Consider $L=f(N)\cup f'(N')\subset M$. It is clear that $L$ is connected. From part(1), For any $z\in L$  there exists a slice $(B(0,\d),\Phi)$ such that $\Phi(0)=z$ so we get a covering of $L$ by slices.  On the other hand, if we have two slices $(B(0,\d),\Phi)$ and $(B(0,\d'),\Phi')$ so that $\Phi(B(0,\d))\cap \Phi'(B(0,\d'))$ is not empty, then in keeping with  part (1), the restriction of $\Phi^{-1}\circ \Phi'$ to ${\Phi'}^{-1}[\Phi(B(0,\d))\cap \Phi'(B(0,\d'))]$ is a diffeomorphim on ${\Phi}^{-1}[\Phi(B(0,\d))\cap \Phi'(B(0,\d'))]$. So we get a structure of connected Banach manifold on $L$, modelled on $\tilde{\D}$. Moreover, by construction, the natural inclusion $i_L:L\ap M$ is injective.
As each slice is an integral manifold of ${\cal D}$ modelled on $\tilde{\D}$,  it follows that,  that $(L,i_L)$ is an integral manifold of $\cal D$.\\
\end{proof}

It remains to show that any equivalent class $L(x)$ of $x\in M$ carries a natural structure of connected Banach manifold modelled on $\tilde{\cal D}_x$. Note that $L(x)$ is the union of all the subset $f(N)$ where $(N,f)$ any integral manifold  through $x$. So $L(x)$ is connected. Moreover, as in the proof of part(2) Lemma \ref {interplaq}, we can cover $L(x)$ by slices and this gives rise to a natural structure of connected Banach manifold on $L(x)$. Again, $(L(x),i_{L(x)})$ is an integral manifold of ${\cal D}$ through $x$, which is maximal by construction.\\

%%%%%%%%%%%%%%%%%%%%%%%%%%%%%%%%%%%%%%%%%%%%%%%%%%%%%%%%%%%%%
\section{Integrability and  Lie invariance}
%%%%%%%%%%%%%%%%%%%%%%%%%%%%%%%%%%%%%%%%%%%%%%%%%%%%%%%%%%%%%%
\subsection{Case of lower trivial weak  distribution}\label{intlieinvlw}
%%%%%%%%%%%%%%%%%%%%%%%%%%%%%%%%%%%%%%%%%%%%%%
{\it In this section we shall adopt the material and arguments used in \cite{St} and \cite{ChSt}}.\\

\noindent Let be $X$ and $Y$ be smooth vector fields on an open set $U$ in a Banach space $E$. Given any integral curve $\g:I=[\a,\b]\ap U$ of $X$ it is well known that the Lie bracket at  some $\g(t)$ is given by
\begin{eqnarray}\label{lieb}
[X,Y](\g(t))=\dis\frac{d}{dt}Y(\g(t))-DX(Y(\g(t))
\end{eqnarray}
where $DX$ is the differential of $X$. Note that, for any diffeomorphism $\phi:U\ap V$, according to the "chain rule" in differentiation, the same type of formula for $\bar{X}=\phi_*X$, $\bar{Y}=\phi_*(Y)$ and the integral curve $\bar{\g}=\phi\circ\g$ of $\bar{X}$ is compatible with as $D\phi([X,Y])(\g(t))=[D\phi\circ X,D\phi\circ Y](\phi\circ\g(t)$ . On the other hand, (\ref{lieb}), depends only of the map of $Y\circ \g$. It follows that the following definition does  not depend of the choice of the chart:

\begin{defi}\label{liebracket}
Let be $X$ an vector field on $M$ and $\g:I=[\a,\b]\ap U\subset M$ an integral curve of $X$ whose range is contained in a chart domain $U$ of a chart $(U,\phi)$. Given any vector field $Y$ along $\g$ (i. e. $Y:I\ap TM$ such that $Y(t)\in T_{\g(t)}M$) the Lie bracket $[X,Y]$ along $\g$ is the vector field characterized by
$$\phi_*[X,Y](\g(t))=\dis\frac{d}{dt}\phi_*Y(\g(t))-D\phi_*X(Y(\g(t))$$
\end{defi}

Note that the previous definition do not depends of the choice of the chart. So  given any $X\in {\cal X}(M)$ and any integral curve $\g:I\ap $Dom$(X)$, the Lie bracket $[X,Y]$ is well defined along $\g$, for any vector field $Y$ along $\g$.

On the other hand, given any smooth curve $\g:I=[\a,\b]\ap M$, we  denote  by $T_\g M$ the restriction of $TM$ to $\g(I)$. For any Banach space $G$, we denote by $L_\g (G,TM)$ the bundle, over $\g$, of morphisms  between the trivial bundle $G\times I$ and $T_\g M$ . \\

%\noindent
 Let ${\cal D}$ be a lower trivial weak  distribution on $M$. Consider a local vector field $X$  and an integral curve  $\g : [\a,\b]\ap M$  of $X$.

\begin{defi}\label{Liederiv}${}$
\begin{enumerate}
\item An upper trivialization of $\cal D$ over $\g$  is   a smooth map $\psi:[\a,\b]\ap L_\g(G,TM)$ (where  $G$ is some Banach space) such that, for each $t\in[\a,\b]$,   the corresponding morphism $\psi_t\in L(G, T_{\g(t)}M)$ satisfies $\psi_t(G)={\cal D}_{\g(t)}$.
\item  Given an upper trivialization $\psi$ as in part 1 , for each $v\in G$, denote by $\psi[v]$ the vector field along $\g$ defined  $\psi[v](t)=\psi_t[v]$. The Lie derivative of $\psi$ by $X$ along $\g$ is defined by:
\begin{eqnarray}\label{lie}
(L_{X}\psi)_t(v)=[X,\psi(v)](\g(t))
\end{eqnarray}
\end{enumerate}
\end{defi}

Remark that, with the previous notations, the map $v\ap (L_{X}\psi)_t(v)$ is linear so we get a smooth map $L_X\psi:I\ap L_\g(G,TM)$ so that for $t\in I$,  $(L_X\psi)_t\in L(G,T_{\g(t)}M)$.

%%\begin{rem}\label{lie prop}${}$
%Definition \ref{Liederiv} is independent of the choice of the chart and so (\ref{lie}) can be defined along any integral curve not necessary contained in a chart domain.
%\end{rem}

\begin{defi}\label{lieinv}${}$\\
Let ${\cal D}$ be a lower trivial weak distribution.
\begin{enumerate}
\item Let $\Theta:\tilde{\cal D}_x\times V\ap TM $, be a lower trivialization around $x$, and $X_u=\Theta(u,\;)$ a lower section on $V$.  \\ We say that the weak  distribution ${\cal D}$ is {\bf Lie invariant by $X_u$}, if , for any $y\in V$, we can find  $\varepsilon>0$, so that, for all $0<\t<\varepsilon$, there exists:\\
$\bullet$ an  upper trivialization $\psi:[-\t,\t] \ap  L_\g(G,TM)$  of $\cal D$   over $\g (t)=\phi^{X_u}_t(y)$ for $t\in [-\t,\t]$,\\
$\bullet$ a smooth field of operator $\L:[-\t,\t]\ap L(G,G)$ \\
which satisfy
 \begin{eqnarray}\label{condLieinv}
L_{X_u}\psi=\psi\circ \L
\end{eqnarray}
\item The  weak distribution ${\cal D}$ is called {\bf  Lie invariant } if for any $x\in M$ there exists a lower trivialization $\Theta:\tilde{\cal D}_x\times V\ap TM$ such that, for any $u\in {\cal D}_x$, the weak distribution  $\cal D$ is Lie invariant by $X_u=\Theta(u,\;)$.
\end{enumerate}
\end{defi}

As in \cite{ChSt}, we have the following Theorem but without the assumption of closeness and existence of a complement for all subspaces ${\cal D}_x$

\begin{theor}\label{thII}${}$\\
Let ${\cal D}$ be a lower trivial  weak distribution. The following properties are equivalent:
\begin{enumerate}
\item ${\cal D}$ is integrable;
\item ${\cal D}$ is Lie invariant;
\item ${\cal D}$ is $ {\cal X}^-_{\cal D}-$invariant.
\end{enumerate}
\end{theor}

\begin{obs}\label{balan}${}$
\begin{enumerate}
\item  On  finite dimensional manifold $M$ consider a "smooth" or "differential" distribution $\cal D$ (see Observations \ref{smoothG} part 3). In a local context  the follwing version of Theorem 4.7 of \cite{Ba} (which is its central result) can be reformulated  in the following way  (Ballan's  proof leads exactly  to this  result). 

\begin{theo} \label{bal}\cite{Ba}${}$\\
we have  equivalence between integrability of $\cal D$ and the following assumptions:

for any $X$ tangent to $\cal D$ and any $x\in$Dom$(X)$, there exists an open neighbourhood $V\subset$Dom$(X)$ and $\varepsilon>0$, a finite set $\{X_1,\cdots,X_p\}$ defined on $V$ and tangent to $\cal D$  and smooth functions $\l_{ij}:]-\varepsilon,\varepsilon[\ap \R$ ($1\leq i,j\leq p$) such that
  \begin{enumerate}
  \item[(a)] $X_1(x),\cdots, X_p(x)$ span ${\cal D}_x$
  \item[(b)] for any $t\in ]-\varepsilon,\varepsilon[$, and $1\leq j\leq p$, $[X,X_i](\g(t))=\dis\sum_{j=1}^p\l_{ij}(t)X_j(\g(t))$ on $ ]-\varepsilon,\varepsilon[$, where $\g(t)=\phi^X_t(x)$
  \item[(c)] ${\cal D}_{\g(t)}$ is generated by $X_1(\g(t)),\cdots, X_p(\g(t))$.
  \end{enumerate}
  \end{theo}
  \smallskip
\textsl{We will show that,  under the previous assumptions (a), (b), (c),  the condition 2 of Theorem \ref{thII} is fulfilled and so, we get that $\cal D$ is integrable.\\
Note that using Proposition \ref{lem1}, the property of lower triviality and  the stability by Lie bracket on an integral manifold we obtain easily the inverse implication.}\\

Indeed,  assume that  previous assumptions (a), (b), (c),  are  fulfilled. At first consider a lower trivialization $\Theta: {\cal D}_x\times V\ap TM$ associated to a  family of linear independent vector fields $\{Z_1,\cdots, Z_n\}$ on a neighbourhood $V$ of $x$ in $M$ (see Observations \ref{smoothG} part 3). Fix some lower section  $X_u=\dis\sum_{j=1}^n u_j Z_j$ for some $u=\dis\sum_{j=1}^n u_j Z_j(x)\in {\cal D}_x$ and consider the integral curve $\g(t)=\phi^{X_u}_t(x)$ on an interval $]-\varepsilon,\varepsilon[$ as in the previous assumptions. Denote by $G$ the vector space $\R^p$ and $\{e_1,\cdots,e_p\}$ its canonical basis. Choose some $0<\t<\varepsilon$. On $[-\t,\t]$ , we consider  the smooth field of linear operators : $t\ap\psi_t$,  from $G$ to $T_{\g(t)}M$, defined by:

$\psi_t(v)=\dis\sum_{j=1}^p v_j X_j(\g(t))$ if $v=\dis\sum_{j=1}^p v_j X_j(x)\in {\cal D}_x$.

 From assumption (c) and Observations \ref{smoothG} part 3, we get an upper trivialization of $\cal D$ along $\g$ .\\ From assumption (b), we can write

$[X,X_i](\g(t))=\dis\sum_{j=1}^p\l_{ij}(t)X_j(\g(t))$ on $ [-\t,\t]$

Consider the field $t\ap \L_t$ of endomorphisms of $G$ defined by 

$\L_t(v)=\dis\sum_{l,j=1}^p\l_{ij}(t)v_l e_j$ for any $v=\dis\sum_{l=1}^pv_l e_l$.

But we have 
 
 $(L_{X_u}\psi)_t(v)=[X_u,\psi(v)](\g(t))=\dis\sum_{l,j=1}^p\l_{ij}(t)v_lX_j(\g(t))=\psi_t\circ\L_t(v)$

So, we get the Lie invariance of $\cal D$ and we can apply Theorem \ref{thII}.

${}$\hfill$\D$
 \item In the same context of  finite dimensional manifolds, recall the remark of Balan about the proof of Theorem 4.1 of \cite{Su}. Consider  the two following  conditions :
 
 "(e)" for any $x\in M$ there exists vector fields $X_1,\cdots X_p$ defined on some neighbourhood of $V$ of $x$ such that
 
 (1) ${\cal D}_x$ is generated by $X_1(x),\cdots X_p(x)$
 
 (2) for every vector field $X$ defined on $V$, there exists $\varepsilon>0$  and smooth functions $\l_{ij}:]-\varepsilon,\varepsilon[\ap \R$ ($1\leq i,j\leq p$) such that $[X,X_i](\g(t))=\dis\sum_{j=1}^p\l_{ij}(t)X_j(\g(t))$ on $ ]-\varepsilon,\varepsilon[$, where $\g(t)=\phi^X_t(x)$
 
 "(f)"  if $\cal D$ is generated by some set $\cal X$ then $\cal D$ is ${\cal X}$-invariant.

Balan points out that the  implication "(e) implies  (f)"   is wrong and, of course, the implication  "(e) implies the integrability of $\cal D$" is wrong too.

The condition (e) can be be replaced by the previous assumptions (a), (b) and (c) proposed by Balan in  \cite{Ba} (see   part 1 of this obervation) which gives rise to a correct version of this Theorem.\\
 In fact  in \cite{Ba}, Balan proposes to replace condition (e) by the following condition 

 (e') for any $x\in$Dom$(X)$, there exists an open neighbourhood $V\subset$Dom$(X)$ and $\varepsilon>0$, a finite set $\{X_1,\cdots,X_k\}$ defined on $V$ and tangent to $\cal D$ such that   
  \begin{enumerate}
  \item[(a)] $X_1(x),\cdots, X_k(x)$ span ${\cal D}_x$
  \item[(b)]  for any vector field $X$ defined on $V$ and tangent to $\cal D$, there exists  smooth functions $\l_{ij}:]-\varepsilon_X,\varepsilon_X[\ap \R$ ($1\leq i,j\leq p$) such that for any $t\in ]-\varepsilon_X,\varepsilon_X[$, and $1\leq j\leq p$, $[X,X_i](\g(t))=\dis\sum_{j=1}^p\l_{ij}(t)X_j(\g(t))$ where $\g(t)=\phi^X_(x)$
  and where\\ $\varepsilon _X=\sup\{\d\;:\; \d\leq\varepsilon,\textrm{ and } \phi^X_t(x)\in V \;\forall |t|<\d\}$
  \end{enumerate}
  
  The essential difference between (e) and (e') is that the parameter  "$\varepsilon $" in  (e') is "maximal"  while in (e), it is not the case.
The corrected proof proposed by Balan, for  the implication   "(e')  implies the integrability of $\cal D$",  contains   implicitly the fact that (e') implies more or less the previous assumptions (a), (b) and (c).
\item Given a set $\cal X$ of vector fields on a finite dimensional manifold $M$, in \cite{St},  Stefan defines the locally subintegrability of $\cal X$  at some $x\in M$ if there exists an open neighbourhood $V$ of $x$ in $M$ and a subset ${\cal X}^\flat$ of $\cal X$ such that

(LS1) ${\cal D}_x=({\cal D}^\flat)_x $  and ${\cal D}^\flat$ is integrable on $V$ where $\cal D$ (resp. ${\cal D}^\flat$) is the distribution generated by $\cal X$ (resp.  ${\cal X}^\flat$)

(LS2) for any $X\in {\cal X}$ there exists $\varepsilon>0$ such that
$$T\phi^X_t(({\cal D}^\flat)_x)=({\cal D}^\flat)_{\phi^X_t(x)}$$
for any $t\in]-\varepsilon,\varepsilon[$.

On the other hand, denote by ${\cal X}^\sharp$ the set of finite  linear combinations of type $\sum\l_iX_i$ where $\l_i$ are smooth functions and $X_i$ belongs to $\cal X$. Then  the Theorem 4 of \cite{St} says:\\
 ${}\;\;\;\;\;\;\;\;\cal D$ is integrable if and only if ${\cal X}^\sharp$ is sub integrable

 According to Balan's remark,  the condition (LS2) is not sufficient , again because "$\varepsilon>0$" given in (LS2) depends of $X$. Moreover, Balan gives a contre example to this Theorem. 
 \end{enumerate}
\end{obs}

\noindent\begin{proof}\so{\it Proof of Theorem 2}. ${}$\\

\noindent According to Theorem 1,  we have only to prove the equivalence $(2) \Longleftrightarrow (3)$.\\

\noindent Assume that  ${\cal D}$ is $ {\cal X}^-_{\cal D}$-invariant. Let $x\in M$ be a fixed point and choose a lower trivialization $\Theta: \tilde{{\cal D}}_x\times V\ap TM$.  Consider a lower section $X_u=\Theta(u,\;)$ and $y\in V$. Note that  there exists $\varepsilon >0$ such that the integral curve  $t\mapsto\phi^{X_u}_t(y)$  of $X_u$ through $y$ is defined for all $t\in ]-\varepsilon,\varepsilon[$.  Choose any $0<\t<\varepsilon$ and  set $\g(t)=\phi^{X_u}_t(y)$ for $t\in [-\t,\t]$. From our assumption, we have $T_y\phi^{X_u}_t({\cal D}_{y})={\cal D}_{\g(t)}$. If $i_y:\tilde{\cal D}_y\ap {\cal D}_y$ is the natural inclusion,  denote by $\psi_t=T_y\phi^{X_u}_t\circ i_y$.  Set $\G_t= T_y\phi^{X_u}_t$. It is clear that $\psi$  is an upper trivialization of ${\cal D}$ over $\g$ with $G=\tilde{\cal D}_y$. On the other hand, we have:

$$(L_{X_u}\psi)_t= [\dot{\G}_t- DX_u(\g(t))\circ \G_t]\circ i_y$$

\noindent But,  we have $\dot{\G}_t= DX_u(\g(t))\circ \G_t$  (see proof of Lemma \ref{invariance}).  So we obtain   $L_{X_u}\psi=0$ on $[-\t,\t]$. Taking $\L=0$ in (\ref{condLieinv}) we get Lie invariance for $X_u$\\

For the converse, as in \cite{St} \cite{ ChSt} and \cite{Nu}, we need the following result whose proof is somewhat different (each space ${\cal D}_x$ can be not closed here)

\begin{lem}\label{localinv}${}$\\
Let $X$ be a local vector field and $\psi$ an upper trivialization of $\cal D$ defined over an integral curve $\g:]-\varepsilon,\varepsilon[\ap V=$dom$(X)$. Moreover we assume that, for any $0<\t<\varepsilon$ there exists a smooth field $\L:[-\t ,\t]\ap L(G,G)$ such that
 $$L_{X}\psi =\psi\circ \L.$$
Then, we have  $T_y\phi^{X}_t[{\cal D}_y]={\cal D}_{\phi^{X}_t(y)}$ for all $0<|t|<\varepsilon$
\end{lem}

Now, assume that ${\cal D}$ is Lie invariant by $X_u$; let us fix some $y\in V=$ Dom$(X_u)$, and  take a maximal integral curve $\g:]\a_y,\b_y[\ap V$ of $X_u$. Consider the set $I=\{t\in ]\a_y,\beta_y[ : T_y\phi^{X_u}_t[{\cal D}_y]={\cal D}_{\phi^{X_u}_t(y)}\}$. This set is clearly  open according to Lemma \ref{localinv}. 
%So, we have $I=]\a_y,\b_y[$. 
Take a sequence $(t_n)$ in $I$ which converges to some $t\in ]\a_y,\b_y[$. From the assumption (2) of the Theorem,  and Lemma \ref{localinv} applied  at the point $\phi^{X_u}_t(y)$, we have $\phi^{X_u}_s({\cal D}_{\phi^{X_u}_t(y)})={\cal D}_{\phi^{X}_{t+s}(y)}$ for $s$ in some neighbourhood $]t-\eta,t+\eta[$ of $t$. As we have some $t_n$ which belongs to $]t-\eta,t+\eta[$, we get that $I$ is closed. So $I= ]\a_y,\b_y[$ and finally we deduce that ${\cal D}$ is invariant by $X_u$.\\
\end{proof}

\noindent \begin{proof}\so{\it Proof of Lemma \ref{localinv}}${}$\\

\noindent  Let $\psi:[-\t,\t]\ap L_\g(G,TM)$ be an upper trivialization of $\cal D$ over an integral curve $\g$ of $X$ such that $\g(0)=y\in V$. Consider any smooth field of operators $\s:[-\t,\t]\ap L(G,G)$ and set $\tilde{\psi}=\psi\circ \s$. On a chart domain, we have
\begin{eqnarray}\label{derivliecompo}
L_{X}\tilde{\psi}=\dot{\tilde{\psi}}-DX\circ\tilde{\psi}=\dot{\psi}\circ \s+\psi\circ\dot{ \s}-DX\circ\psi\circ \s=L_{X}\psi\circ \s+\psi\circ\dot{\s}
\end{eqnarray}

\noindent Assume that $L_{X}\psi=\psi\circ \L$ for some smooth field of operators $\L:[-\t,\t]\ap L(G,G)$. Then we have:
\begin{eqnarray}\label{choixs}
L_{X}\tilde{\psi}=\psi\circ\L\circ \s+\psi\circ \dot{\s}=\psi\circ[\L\circ \s +\dot{\s}]
\end{eqnarray}

\noindent Consider the  solution ({\it again denoted  by $\s$})  of the linear equation $\dot{\s}=(-\L)\circ \s$  with initial condition $\s_0=Id_G$.  So $\s$ is a smooth field of isomorphisms of $G$ and in particular, for this choice of $\s$, we have $\tilde{\psi}(t)[G]={\cal D}_{\g(t)}$ for any $t\in [-\t,\t]$. Moreover, using (\ref{choixs}), we have $L_{X}\tilde{\psi}=0$.  
In fact the relation (\ref{choixs}) do not depends of the choice of the chart, so we can get the same result even if $\g([-\t,\t])$ is not contained in a chart domain.\\
So  we can assume that there exists an upper trivialization $\psi :[-\t,\t]\ap L(G,TM)$ such that $L_{X_u}\psi=0$ on $\g$. Again we set  $\G_t= T_y\phi^{X}_t$. Then  $\S_t=[\G_t]^{-1}\circ\psi_t$ is a smooth field of continuous operators from $G$ to $E\equiv T_xM$ defined on $[-\t,\t]$.

\noindent On a chart domain we have

$$\dot{\psi}=\dot{\G}\circ \S+\G\circ\dot{\S}= DX_u(\g)\circ\G\circ\G^{-1}\circ \psi+\G\circ\dot{\S}=DX_u\circ \psi+\G\circ\dot{\S}$$

%\noindent It follows that, on $[-\t,\t]$ we have

%$\dot{\psi}=\dot{\G}\circ \S+\G\circ\dot{\S}= DX_u(\g)\circ\G\circ \psi+\G\circ\dot{\S}=DX_u\circ \psi+\G\circ\dot{\S}$

\noindent According to (\ref{lieb}) and  (\ref{lie}), on $[-\t,\t]$ we have
\begin{eqnarray}\label{LXpsi}
L_{X}\psi=\G\circ\dot{\S}
\end{eqnarray}

Now, (\ref{LXpsi}) do not depend of the choice of the  chart. As $\G$  and $\dot{\S}$ are  intrinsically defined,  so even if $\g([-\t,\t])$ is not contained in a chart domain, we can  obtain the same relation.\\
 Now from our assumption  $L_{X}\psi=0$, as $\G_t$ is an isomorphism, we must have $\S_t=\S_0=\psi_0$. We conclude that, for any $t\in [-\t,\t]$, we have  $[\G_t]^{-1}\circ \psi_t(G)=\psi_0(G)={\cal D}_y$ and finally
\begin{eqnarray}\label{invtau}
T_y\phi^{X}_t[{\cal D}_y]=\psi_t(G)={\cal D}_{\g(t)}.
\end{eqnarray}
Now, according to the assumption in this Lemma, there exists $\varepsilon >0$ such that, we are in the previous situation for any interval $[-\t,\t]$ with $0<\t<\varepsilon$.  So  (\ref{invtau}) is true for any $0<|t|<\varepsilon$
\end{proof}

%%%%%%%%%%%%%%%%%%%%%%%%%%%%%%%%%%%%%%%%%%%%%%%%%%%%%%%%%%%%%%
\subsection{Case of  upper trivial weak distribution}\label{intlieinvstup}
%%%%%%%%%%%%%%%%%%%%%%%%%%%%%%%%%%%%%%%%%%%%%%

Let ${\cal D}$ be an upper trivial weak  distribution on $M$ (see subsection \ref{preliminaires et resultats}). By analogy with lower sections (see subsection \ref{resultats}),  for any  upper trivialization  $\Psi :F\times  V\ap TM$  such that the associated lower trivialization $\Theta : \tilde{{\cal D}}_x\times V\ap TM$, an {\bf upper section} is a local vector field on $M$ defined by
  \begin{eqnarray}\label{uppersec}
Z(y)=\Psi(u,y)  \textrm{ for any } u\in F
\end{eqnarray}

\begin{rem}\label{generateur}${}$
\begin{enumerate}
\item The set ${\cal X}^+_{\cal D}$ of upper sections generates ${\cal D}$. 
\item Given any  upper trivialization $\Psi:F\times V\ap TM$ at $x$, consider the module ${\cal X}_{\cal D}(V)$ of vector fields $X\in {\cal X}(M)$ whose domain contains $V$ and which are tangent to $\cal D$ on $V$. The set ${\cal X}^+_{\cal D}(V)=\{Z_v=\Psi(v,.),\; v\in F$ is contained in ${\cal X}_{\cal D}(V)$ and of course $\{Z_v(y),\;v\in F\}$ generates ${\cal D}_y$ for all $y\in V$.  Moreover, if $F$ has a Schauder basis $\{e_\a,\a\in A\}$, then  the convex hull of $\{Z_{e_\a}(y),\a\in A\}$  is  dense in $\tilde{\cal D}_y$
\item If $\Theta:\tilde{\cal D}_x\times TM$ is the lower trivialization associated to the  upper trivialization $\Psi$ then each lower section $X_u=\Theta(u,.)$ can be written $X_u=\Theta(\Psi(u,x),\,)$ and so the set ${\cal X}^-_{\cal D}(V)$ of such lower sections is contained in  ${\cal X}^+_{\cal D}(V)$. 
\end{enumerate}
\end{rem}

Let  $\cal D$ be a  upper trivial  weak distribution on $M$.
Let $V$ be the domain of a chart around $x\in M$. Consider a  upper trivialization $\Psi:F\times V\ap TM$ and $\Theta:\tilde{{\cal D}}_x\times V\ap TM$ the associated lower section.  Given any smooth function $\s:V\ap F$, let $Z_\s=\Psi(\s,\;)$ be the associated vector field on $V$. Consider $\g:[-\t,\t]\ap V$  an integral curve of $Z_\s$, then, $\Psi_\g$, defined by $\Psi_\g(t)[v]=\Psi{(v,\g(t))} $, is an upper trivialization of $\cal D$ along $\g$,  according to Definition \ref{Liederiv}. So the Lie derivative of $\Psi$ by $Z_\s$ along $\g$  is $L_{Z_\s}\Psi_{\g}$ which we simply denoted by $L_{Z_\s}\Psi$. We can also define directly $L_{Z_\s}\Psi$ by:

%\begin{rem}\label{liestrong}${}$
%\begin{enumerate}
%\item As in subsection \ref{intlieinvlw},  the definition of  $L_{Z_\s}\Psi$ is independent  of  the choice of the chart and can be defined along any integral curve of $X_u$.\\
%\item
% If $\Psi:F\times V\ap TM$ is an upper trivialization   we have  (see \cite{St}):
\begin{eqnarray}\label{strongliebrac}
L_{Z_\s}\Psi(v,\g(t))=[Z_\s,Z_v](\g(t)) \textrm{ for any } Z_v=\Psi(v,\;)
\end{eqnarray}
%\end{enumerate}
%\end{rem}

\begin{defi}\label{stlieinv}${}$\\
An  upper trivial weak distribution ${\cal D}$ is called Lie bracket invariant  if, for any $x\in M$, there exists an upper trivialization $\Psi:F\times V\ap TM$ such that for any $u\in F$ , there exists $\varepsilon>0$,  such that, for all $0<\t<\varepsilon$, there exists
 a smooth field of operators $\L:[-\t,\t]\ap L(F,F)$ with the following property
 \begin{eqnarray}\label{condLieinvst}
L_{X_u}\Psi=\Psi\circ \L
\end{eqnarray}
along the integral curve $t\mapsto \phi^{X_u}_t(x)$ on $[-\t,\t]$ of any lower section $X_u=\Theta(\Psi(u,x),\;)$.
\end{defi}

\begin{rem}\label{Liebracket}${}$
%\begin{enumerate}\item
According to (\ref{strongliebrac}),  the property (\ref{condLieinvst}) is equivalent to
\begin{eqnarray}\label{condition crochet}
 [X_u,Z_v](\g(t))=\Psi(\L_t(v),\g(t)) \textrm{ for any } Z_v=\Psi(v,\;)
\end{eqnarray}
along $\g(t)=\phi^{X_u}_t(x)$.\\
 (\ref{condition crochet}) justifies the term "Lie bracket invariant" in Definition \ref{stlieinv}.
%\item In fact, in Definition \ref{stlieinv}, we can impose that the morphism $\L(t)$ in (\ref{condLieinvst}) satisfies
%$$ \ker \L(t)=\ker\Psi_{\g(t)}$$
% (see the proof of Theorem \ref{thIII}).
%\end{enumerate}
\end{rem}
With these definitions we have:

\begin{theor}\label{thIII}${}$\\
Let ${\cal D}$ be an  upper trivial weak  distribution. The following propositions are equivalent:
\begin{enumerate}
\item ${\cal D}$ is integrable;
\item ${\cal D}$ is Lie bracket invariant;
\item ${\cal D}$ is $ {\cal X}^-_{\cal D}$-invariant.
\end{enumerate}
\end{theor}
\begin{rem}\label{B}${}$
\begin{enumerate}
\item The assumption "the kernel of $\Psi$ is complemented in each fiber" is automatically satisfied if the kernel of $\Psi$ is finite dimensional or finite codimensional in each fiber, or in the context of Hilbert manifold.  
\item Recall that when $M$ is a finite dimensional manifold, if, for any $x\in M$,  each module ${\cal X}_{{\cal D}_x}$   of germs of vector fields  is finitely generated then ${\cal D}$ is  upper trivial(see Observations \ref{smoothG} part 3.) 
So, the   formulation Theorem 4.7 of \cite{Ba} can be given in its original way:

  if  each  ${\cal X}_{{\cal D}_x}$ is finitely generated for any $x\in M$,  then $\cal D$ is integrable if and only if $\cal D$ satisfies  the properties (a), (b) and (c) in Theorem \ref{bal}. \\
  
  Of course, this result is a direct consequence of this last Theorem, but, we can easily  see that this Theorem 4.7 can be directly deduced from Theorem \ref{thIII}. This proof is left to the reader.
\end{enumerate}
\end{rem}

Coming back to the context of Corollary \ref{fibre},  let  $\Pi:{\cal F}\ap M$ be a Banach fiber bundle over $M$ with  typical fiber $F$, $\Psi:{\cal F}\ap TM$ a morphism of bundle whose kernel is complemented in each fiber. We denote by ${\cal S}({\cal F})$ the set of local sections of  $\Pi:{\cal F}\ap M$ , that is smooth maps $\s:U\subset M\ap {\cal F}$ such that $\Pi\circ \s=Id_U$ where $U$ is an open set of $M$. The maximal such open set is called the {\it domain} of $\s$ and denoted    Dom$(\s)$. \\ 

A subset ${\cal S}$ of  ${\cal S}({\cal F})$ is called a {\bf generating upper set } of $\cal D$ if for  any $x\in M$, the set  ${\cal X}_{\cal S}=\{\Psi\circ \s, \; \s\in {\cal S}\}$ contains ${\cal X}^+_{\cal D}$.  Of course ${\cal S}({\cal F})$ is a maximal generating upper set.  We introduce some condition on  ${\cal X}_{\cal S}$ which will be used in the next theorem:

{\it ${\cal X}_{\cal S}$ satisfies the condition (LB)} if: 

\noindent for any local section $\s\in{\cal S}$,  there exists an open set $V\subset  $Dom$(\s)$ on which ${\cal F}$  is trivializable   
and  for  any $x\in V$  we have the following property:

 given any  integral curve $\g:]-\varepsilon,\varepsilon[\ap V$  of $X=\Psi\circ \s$  with  $\g(0)=x$, there exists  a smooth field $\L:]-\varepsilon,\varepsilon[\ap L({\cal F}_x,{\cal F}_x)$ such that 
 \begin{eqnarray}\label{LBF}
[\Psi\circ\s,\Psi(u,\;)](\g(t))=\Psi(\L_t(u),\g(t)) \textrm{ for any } t\in ]-\varepsilon,\varepsilon[ \textrm{ for any }  u\in {\cal F}_x
\end{eqnarray}

Then, using  Theorem 3 we get the following Theorem

\begin{theor}\label{thIV}${}$\\
Let $\Pi:{\cal F}\ap M$ be a Banach fiber bundle over $M$ with  typical fiber $F$ and   $\Psi:{\cal F}\ap TM$ a morphism of bundles such that the kernel of $\Psi$ is complemented in each fiber and denote ${\cal D}=\text{Im }\Psi$.
\begin{enumerate}
\item Then ${\cal D}$ is an  integrable distribution if and only there exists a generating upper set ${\cal S}$  such that ${\cal X}_{\cal S}$ satisfies the condition (LB)

Moreover, when (LB) is  satisfied,  if $S_x$ fulfills ${\cal F}_x=\ker \Psi_x \oplus S_x$, there exists \\ $\L:]-\varepsilon,\varepsilon[\ap L({\cal F}_x,{S}_x)$  which satisfies (\ref{LBF})

\item If ${\cal D}$ is a closed distribution, then this distribution is integrable if and only if (LB) is satisfied where (\ref{LBF}) can be replaced by
\begin{eqnarray}\label{LBFC}
[\Psi\circ\s,\Psi(u,\;)](\g(t))\in \Psi_{\g(t)}(S_x) \textrm{ for any }  t\in ]-\varepsilon,\varepsilon[ \textrm{ for any }  u\in {\cal F}_x
\end{eqnarray}
\end{enumerate}
\end{theor}

\subsection{Proofs of Theorem 3 and Theorem 4}
%%%%%%%%%%%%%%%%%%%%%%%%%%%%%%%%%%%%%%%%%%%%%%%%%%%%%%

\noindent\begin{proof}\so{\it Proof of Theorem 3}${}$\\

According to Theorem 1, we have only to prove $(1)\Longleftrightarrow (2)$.\\

 From  Lemma \ref{localinv}, property 2 of Theorem \ref{thIII}  implies that for any $x\in M$,  we have $T_x\phi^{X_u}_t({\cal D}_x)={\cal D}_{\phi^{X_u}_t(x)}$ for all $t$ such that $(x,t)\in \O_{X_u}$.  From Proposition \ref{immersion}, $(B(0,\d),\Phi)$ is an integral manifold through $x$.  So $(2)\Longrightarrow (1)$.  \\

 For the converse, we will use the following Lemma:

 \begin{lem}\label{compLie}${}$\\
Assume that $\cal D$ is integrable.  Let $\Psi:F\times V\ap TM$ be a  upper trivialization, and $\s:V\ap F$ a smooth map and let $X=\Psi(\s,\;)$ be the associated vector field on $V$. Consider an integral curve $\g:]-\varepsilon,\varepsilon[\ap V$ of $X$ such that $\g(0)=x$. Then there exists a smooth field $\L:]-\varepsilon,\varepsilon[ \ap L(F,F)$ such that :
 $$L_X\Psi(v,\g(t))=\Psi(\L_t(v),\g(t))$$
 \end{lem}

\noindent So, for $\s(y)=(u,y)$, with $u\in S$, the vector field $Z_\s$ is the lower section $X_u$ for $u\in S$ et clearly Lemma \ref{compLie} proves  $(1)\Longrightarrow (2)$.\\
\end{proof}\\
 \noindent\begin{proof}\so{\it Proof of Lemma \ref{compLie}}${}$\\

Recall that we have assumed  that ${\cal D}$ is integrable. Fix some $x\in M$.  Take an  upper trivialization $\Psi:F\times V \ap TM$  around $x$  and let $\Theta:\tilde{\cal D}_x\times V \ap TM$ be the associated lower trivialization. We can choose $V$ such that $TM_{|V}\equiv E\times V$. Recall that we have a decomposition $F=\ker \Psi_x\oplus S$, and   ${\Theta}=({\theta}_y\circ [{\theta}_x]^{-1},\; )$ where $\theta_y$ is the restriction to $S$ of $\Psi_y$ (see  the proof of Proposition \ref{upperloc}). At first note that  any lower section $X_u$ can be written $X_u=\Theta(\Psi(u,x),\;)$ for some $u\in F$  and also $X_u=\theta(u,\;)$ but  with $u\in S$. On the other hand,  according to  Lemma \ref{interplaq},  $(B(0,\d),\Phi)$  is an integral manifold of ${\cal D}$ through $x$ (associated to $\Theta$). In particular, $\tilde{\theta}_y$ is an isomorphism from $S$ to $\tilde{\cal D}_y$. Given $u\in F$, there exists an unique $v\in S$ such that $\Psi_y(u)=\theta_y(v)$ so $u\in \ker \Psi_y\oplus S$. It follows that $F=\ker \Psi_y\oplus S$.

Set $N=\Phi(B(0,\d))\subset M$ endowed with its Banach manifold structure. Without loss of generality, we can identify $S$ with $\tilde{\theta}_x(S)=\tilde{\cal D}_x$ and so we can consider   $N$ as an open set of $\tilde{\cal D}_x$.  Denote by $i:\tilde{\cal D}_x\ap T_xM\equiv E$ the canonical inclusion.
We have $T_yN \equiv S\times\{y\}$.  On $N$,  each $\tilde{\Psi}_y$ can be considered as an  element of $L(F,S)$. By arguments similar to those used in the proof of Lemma \ref{Xdiff}, we can show that  $y\mapsto\tilde{\Psi}_y$   is a smooth field of operators in $L(F,S)$. So, $y\mapsto \tilde{\theta}_y$ is  a smooth field of isomorphisms of $S$. We get a smooth  field $\pi_y=[\tilde{\theta}_y]^{-1}\circ \tilde{\Psi}_y$ of operators in $L(F,S)$ with the following properties:
\begin{eqnarray}
\label{decompproj} \Psi_y=\theta_y\circ \pi_y\\
\label{kerproj} \ker \pi_y=\ker\tilde{\Psi}_y=\ker \Psi_y\\
\label{proj} \pi_y(u)=u \textrm{ for all } u\in S
\end{eqnarray}
Take any smooth map $\s:V\ap F$. Then $ Z_\s(y)=\Psi_y\circ\s(y)$ for $y\in V$ (resp. $\tilde{ Z}_\s(y)=\tilde{\Psi}_y\circ\s(y)$ for $y\in N$) is a smooth vector field on $V$ (resp.  on $N$) and we have the relations:
\begin{eqnarray}\label{upperlower}
 \Psi(\s(y),y)=Z_\s(i(y))=i[\tilde{Z}_\s(y)]=i\circ \tilde{\theta}_{i(y)}\circ\pi_{y}(\s(y))=\theta_{i(y)}\circ \pi_{y}(\s(y))=\theta(\pi_y(\s(y)),y)
\end{eqnarray}

 Consider the integral curves $\g(t)=\phi^{Z_\s}_t(x)$ and $\tilde{\g}(t)=\phi^{\tilde{Z}_\s}_t(x)$ for $t\in ]-\varepsilon,\varepsilon[$. Of  course we have $\g(t)=i\circ\tilde{\g}(t)$. For simplicity, we set:

 $\s(\g(t))=\s(t)$ and $\s(\tilde{\g}(t))=\tilde{\s}(t)$

   $P(t)=\pi_{\tilde{\g}(t)}$. \\
Note that , using  (\ref{upperlower}) we have
\begin{eqnarray}\label{factorisationII}
\Psi(v\;,\g(t))=\theta(P(t)[v],\g(t))
\end{eqnarray}

\noindent Now, in keeping with  (\ref{strongliebrac}), for any $v\in S$, we have:
\begin{eqnarray}\label{lietheta}
L_{Z_\s}\Psi(v,\g(t))=[Z_\s, X_v](\g(t))=L_{Z_\s}\theta(v,\g(t)).
\end{eqnarray}
 As we have  $Z_\s=i_*\tilde{Z}_\s$ and $X_v=i_{*}\tilde{X}_v$ , on the Banach manifold  $N$, we get
 $$[\tilde{Z}_\s,\tilde{X}_v](\tilde{\g}(t))\in T_{\tilde{\g}(t)}N\equiv S\times\{\tilde{\g}(t)\}.$$
 Note that we have $[Z_\s,X_v]=i_*[\tilde{Z}_\s,\tilde{X_v}]$. It follows that we have:
 \begin{eqnarray}\label{explietheta}
L_{Z_\s}\theta(v,\g(t))=i_*[\tilde{Z}_\s,\tilde{X_v}](\tilde{\g}(t))
\end{eqnarray}

Without loss of generality, we can choose $\d>0$ small enough such that on $N=\Phi(B(0,\d))$ we have :
  \begin{eqnarray}\label{majoration}
||\tilde{\theta}(\;,y)||\leq K \textrm{ and } ||D_2\tilde{\theta}(\;,y)[\;]||\leq K \textrm{  for all } y\in N.
\end{eqnarray}
 On the other hand, we have:
\begin{eqnarray}\label{liedecomposition}
[\tilde{Z}_\s,\tilde{X}_v](\tilde{\g}(t))=
 D_2\tilde{\theta}(v,\tilde{\g}(t))[\tilde{\theta}(P(t)[\tilde{\s}(t)],\tilde{\g}(t))]-D_2\tilde{\theta}(\tilde{\theta}(P(t)[\tilde{\s}(t)],\tilde{\g}(t))[\tilde{\theta}(v,\tilde{\g}(t))].
\end{eqnarray}

\noindent From (\ref{majoration}), we have :
\begin{eqnarray}\label{majoration1}
|| D_2\tilde{\theta}(v,\tilde{\g}(t))[\tilde{\theta}(P(t)[\tilde{\s}(t)],\tilde{\g}(t))]||\leq K||v||.||\tilde{\theta}(P(t)[\tilde{\s}(t)],\tilde{\g}(t)||\leq K^2||v||.||P(t)||.||\tilde{\s}(t)||\\ \nonumber
\leq K^2.||v||.||P||_{\infty}||\tilde{\s}||_{\infty}.
\end{eqnarray}
 \noindent In the second member of (\ref{liedecomposition}),  the same majoration is true for  $$||D_2\tilde{\theta}(\tilde{\theta}(P(t)[\tilde{\s}(t)],\tilde{\g}(t))[\tilde{\theta}(v,\tilde{\g}(t))||.$$
 \noindent So we get

 $||[\tilde{Z}_\s,\tilde{X}_v](\tilde{\g}(t))|| \leq  2K^2||P||_{\infty}||\tilde{\s}||_{\infty}||v||$\\

\noindent  It follows that , for each $t\in ]-\varepsilon,\varepsilon[$,   the map
 $v\ap [\tilde{Z}_\s,\tilde{X}_v](\tilde{\g}(t))$ is linear continuous map $\tilde{\L}_{\tilde{\g}(t)}$ from $S$ to $S\times\{\tilde{\g}(t)\}$ . We set
 $$\bar{\L}(t)=[\tilde{\theta}_{\tilde{\g}(t)} ]^{-1}\circ \tilde{\L}_{\tilde{\g}(t)}$$
 Clearly, $t\mapsto \overline{\L}(t)$ is a smooth field of endomorphisms of $S$ and taking in account (\ref{explietheta}) we get
 \begin{eqnarray}
L_{Z_\s}\theta(v,\g(t))=\theta(\bar{\L}(t)[v],\g(t))
\end{eqnarray}
Now from (\ref{factorisationII}), with the same argument (\ref{derivliecompo}) used in the  proof of Lemma \ref{localinv},  we get:
\begin{eqnarray}\label{(2)}
L_{Z_\s}\Psi_{\g(\;)}=L_{Z_\s}\{\theta_{\g(\;)}\circ P\}+\theta_{\g(\;)} \circ\dot{P}=\theta_{\g(\;)}\circ P\circ\bar{\L}+\theta_{\g(\;)}\circ \dot{P}
\end{eqnarray}
But, according to the definition of  $P$ and (\ref{proj}), we have $P\circ\dot{P}=\dot{P}$. So using  (\ref{(2)}), we get:
$$L_{Z_\s}\Psi=\Psi(\bar{\L}[\;]+\dot{P}[\;],\;)$$
which ends the proof  of Lemma \ref{compLie} by setting $\L=\bar{\L}+\dot{P}$.\\
\end{proof}

\noindent\begin{proof}\so{\it Proof of Theorem \ref{thIV}}${}$\\

(1) According to the context of the proof of Corollary \ref{fibre}, for any given $x\in M$ we consider  a local trivialization of ${\cal F}$ on an open set $V$ around $x$, so that  the morphism $\Psi$ can be identified,  as a map $\Psi: {\cal F}_x\times V\ap TM$ and ${\cal F}_x=\ker\Psi_x\oplus S_x$ and let $\Theta:S_x\times V\ap TM$ be the associated lower trivialization. In this context, taking any  $\s(y)=(u,x)$, for any $u\in {S}_x$ in (LB)  for any $x\in M$,  we get  property (2) of Theorem 3 so (LB) is a sufficient condition for integrability of ${\cal D}$.\\
Assume now that ${\cal D}$ is integrable and consider  an upper generating set $\cal S$ of  ${\cal D}$  and any section $\s\in {\cal S}$ defined  on an open set $U$. Fix any $x\in U$.   From Lemma \ref{compLie} we get  (\ref{LBF}) with $\L:]-\varepsilon,\varepsilon[\ap L({\cal F}_x,S_x)$ So in particular, if (LB) is true for some $\L:]-\varepsilon,\varepsilon[\ap L({\cal F}_x,{\cal F}_x)$, then $\cal D$ is integrable and by Lemma \ref{compLie} we can find $\L':]-\varepsilon,\varepsilon[\ap L({\cal F}_x,{S}_x)$ which satisfies (\ref{LBF}) .\\

(2) If now ${\cal D}$ is a closed distribution, then in each fiber  ${\cal F}_x=\ker \Psi_x\oplus S_x$, $\theta_x$ is a continuous bijective morphism between both Banach space $S_x$ and ${\cal D}_x$ so  $\theta_x$ is an isomorphism. In particular, $\tilde{\cal D}_x$ and ${\cal D}_x$ are equivalent as Banach spaces. Coming back to the previous local context of the  upper trivialization $\Psi:{\cal F}_x\times V\ap TM$, the map $y \mapsto \theta_y$ is a smooth field of isomorphisms from $S_x$ to $\Psi_{\g(t)}(S_x)$. \\
If ${\cal D}$ is integrable, from (\ref{LBF}) and the properties of $\L$  we obtain (\ref{LBFC}). For the converse, it is sufficient to set
$$\L_t(u)=[\theta_{\g(t)}]^{-1}[Z_\s,Z_u](\g(t))$$
to get (\ref{LBF}).\\
\end{proof}

%%%%%%%%%%%%%%%%%%%%%%%%%%%%%%%%%%%%%%%%%%%%%%%%%%%%%%%%%%%%%%
\section{Applications}\label{appli}
%%%%%%%%%%%%%%%%%%%%%%%%%%%%%%%%%%%%%%%%%%
\subsection{Banach Lie Algebroid}
%%%%%%%%%%%%%%%%%%%%%%%%%%%%%%%%%%%%%%%%%%%
The concept of Lie algebroid was first introduced by J. Pradines  in relation with Lie groupoids (cf  \cite{Pr}). The theory of algebroids was developped  by A. Weinstein (\cite{We}) and, independently, by M. Karasev (\cite{Ka}), in view of the symplectization
of Poisson manifolds. This  theory   has also an important role  as models in mechanic and mathematical physic (for a survey, see \cite{Li} and \cite{Kos2} for instance). On the other hand, this concept of Lie algebroid   can be extend in the infinite dimensional case: in \cite{KisLe} the authors  build   variational Lie algebroids of  the infinite jet bundles over a vector bundle over a finite dimensional manifolds. This construction can be situated in the Frechet manifold framework . In fact, this context is very special : infinite jet bundles over a vector bundle  over finite dimensional  manifolds are projective limit of finite dimensionnal Banach spaces so we get a set of coordinate on such a space. The existence of  these coordinates   is crucial in this construction of the variational Lie algebroid. In this paper, we look for the Banach manifold context and in this framework we  do not have (local)  coordinates in general.  According to the classical definition of a Lie algebroid in finite dimension  we introduce:

\begin{defi}\label{Algebroid}${}$:\\
A Banach Lie algebroid   structure on a Banach bundle $\Pi:{\cal A}\ap M$  is a quadruple ${({\cal A},\Psi,M,\{\;,\;\})}$  such that
\begin{enumerate}
\item    a bracket $\{\;,\;\}$, i. e.  is a composition law $(\s_1,\s_2)\mapsto \{\s_1,\s_2\}$ on the set of global  sections ${ \cal S}(\cal A)$ of  $\Pi: {\cal A} \ap M$, such that,  $({ \cal S}(\cal A),\{\;,\;\})$ has a Lie algebra structure;
\item $\Psi:{\cal A} \ap TM$  is a smooth vector bundle morphism;
 \item the Leibniz property is satisfied:\\
   for any  smooth function  $f$ defined on  $M$ and any sections $\s_1,\s_2\in { \cal S}(\cal A)$ we have :
$$\{\s_1,f\s_2\}=f\{\s_1,\s_2\}+ df(Z_{\s_1})\s_2 $$
where $Z_{\s_1}=\Psi\circ \s_1$  is the vector field associated to $\s_1$.
\end{enumerate}
The quadruplet ${({\cal A}, \Psi,M,\{\;,\;\})}$ is called a Banach Lie algebroid and  $\{\;,\;\}$  (resp. $\Psi$) is called the {\bf Lie bracket} on ${\cal A}$, (resp.  the {\bf anchor morphism}).
\end{defi}

As in finite dimension, the Jacobi identity and the Leibniz property implies that  $\Psi$ gives rise to a  Lie algebra morphism from ${\cal S}({\cal A})$ into ${\cal X}(M)$ i. e.
\begin{eqnarray}\label{alge}
[\Psi\circ \s_1,\Psi\circ\s_2]=\Psi\circ\{\s_1,\s_2\}
\end{eqnarray}

Given some open  set $U$ in $M$, we denote by ${\cal A}_U$ the restriction of the Banach bundle $\Pi:{\cal A}\ap M$ to the Banach manifold $U$: ${\cal A}_U=\Pi^{-1}(U)$; the set of sections of ${\cal A}_U$ will be denote by ${\cal S}({\cal A}_U)$.\\

 In finite dimension, it is classical that a bracket  $\{\;,\;\}$  on a Lie algebroid ${({\cal A},\Psi,M,\{\;,\;\})}$  is compatible the sheaf of sections of $\Pi:{\cal A}\ap M$ or, for short,  is {\bf localizable}  (see for instance \cite{Ma}). By this property, we mean  the following :
\begin{enumerate}
\item[(i)]  for any open set $U$ of $M$, there exists a unique  bracket $\{\;,\;\}_U$ on the space of sections  ${\cal S}({\cal A}_U)$ such that, for any $s_1$ and $s_2$ in  ${\cal S}({\cal A})$, we have:
$$\{{s_1}_{|U},{s_2}_{|U}\}_U=(\{s_1,s_2\})_{| U}$$
\item[(ii)]  (compatibility with restriction) if $V\subset U$ are open sets, then, $\{\;,\;\}_{U}$  induces a unique  bracket  $\{\;,\;\}_{UV}$ on ${\cal S}({\cal A}_V)$ which coincides with  $\{\;,\;\}_V$ (induced by $\{\;,\;\}$).
\end{enumerate}

Using the same arguments as in finite dimension, when $M$ is smooth regular, we can prove that, for  any   Lie algebroid   ${({\cal A},\Psi,M,\{\;,\;\})}$, its bracket is localizable (see \cite{CaPe}). But, if $M$ is not smooth regular, we can no more used this argument. Unfortunately, we have {\it no example of Lie algebroid} for which the Lie bracket is not localizable. Note that, according to \cite{KrMi}  sections 32.1, 32.4,  33.2 and 35.1,  this problem is similar  to the problem  of localization (in an obvious sense) of  global derivations of  the module of smooth functions on $M$  or  the module  of differential forms on $M$. In \cite{KrMi} and, to our known,  more generally in the literature, there exists no example of such derivations which are not localizable. On the other hand, even if $M$ is not regular, the classical Lie bracket of vector fields on $M$ is localizable.  So, there always exists an anchored bundle ${\cal A}=TM$ and a Lie bracket algebroid $(TM, Id,M,[;,.])$ for which its Lie bracket is localizable.  Moreover,  in Examples \ref{alg}  we  do not assume   that $M$ is regular but, nevertheless, these Lie brackets are also localizable.\\  {\bf Thus, in the  Definition \ref{Algebroid},  we moreover impose that, if  $M$ is not regular, then  the Lie bracket  of the Lie  algebroid  is localizable} even if  this assumption could be, in fact, unnecessary when $M$ is not regular.\\

\begin{rem}\label{Kos}${}$:\\
In finite dimension we have many equivalent definitions of a Lie algebroid:  a Lie algebroid structure on a vector bundle ${\cal A}\ap M$ may be characterized by:

$\bullet$ a Lie bracket on an anchored bundle $({\cal A},\r)$;

$\bullet$ a linear Poisson structure on  ${\cal A}^*$;

$\bullet$ a linear Schouten structure on the exterior algebra $\L^\bullet {\cal A}^*$;

$\bullet$ a differential operator $d$ on the module of sections  ${\cal S}(\L^\bullet {\cal A}^*)$ with $d\circ d=0$.

This last approach can  be interpreted  in the context of supermanifolds (see   \cite{Kos1}). It is precisely  this last aspect  which is used in \cite{KisLe} for the construction of variational Lie algebroids.  However, in  the context of Banach manifolds, we have many obstructions in the generalization of  the previous equivalent definitions.  For instance, when the typical fiber of ${\cal  A}$ is an infinite dimensional Banach space,  in general, such a differential operator $d$ could be not  localizable (see \cite{KrMi}  section 35.1).  Moreover,  according to the fact that the bidual $E^{**}$ of a Banach space $E$  may contains strictly $E$,  we must impose complementary conditions on $d$, to get a Lie  algebroid structure on $\cal A$ by this way. On the other hand, the set of sections of  the graded algebra  $\L^{\bullet}{\cal A}$ is not generated by  elements of degree $0$ and $1$. So, we cannot   extend the Lie algebroid bracket to a unique linear  Schouten  bracket  on ${\S}(\L^{\bullet}{\cal A})$,  such that $\{s,f\}=df( \r(s))$ for any  $s\in{\S}({\cal A})$ and any smooth function $f$ on $M$.  Such problems are studied in \cite{CaPe}.\\
\end{rem}

 Now, in the general situation of a Lie Banach algebroid ${({\cal A},\Psi,M,\{\;,\;\})}$, if the kernel of  $\Psi$ is complemented in each fiber ${\cal A}_x$, then the distribution ${\cal D}=\Psi({\cal A})$ is upper trivial. So, from Theorem 4 we then  get:

\begin{theor}\label{thV}${}$\\
  %If a Banach bundle $\Pi:{\cal A}\ap M$ has a structure of  local Banach Lie algebroid associated to an anchor $\Psi:{\cal A}\ap TM$ and if the kernel of $\Psi $ is complemented in each fiber, then ${\cal D}=\Psi({\cal A})$ is an integrable weak distribution.\\
%In particular, 
Let ${({\cal A},\Psi,M,\{\;,\;\})}$ be a  Banach Lie algebroid.  If  the kernel of $\Psi $ is complemented in each fiber, then  ${\cal D}=\Psi({\cal A})$ is an integrable weak distribution.
\end{theor}

 \begin{exe}\label{alg}${}$
 \begin{enumerate}
\item Let $\Pi: {\cal A} \ap M$ be a weak subbubdle of $TM\ap M$  i.e. : $\cal A\subset TM$ , $\Pi$ is the restriction to $\cal A$ of the canonical projection of $TM$ onto $M$, and the canonical   inclusion  $i:{\cal A} \ap TM$ is a morphism bundle. Any  section of  $\Pi: {\cal A} \ap M$ induces a  vector field on $M$. Assume that the set  of sections ${\cal S}({\cal A})$ is stable by Lie bracket of vector fields,  which means that the associated weak  distribution ${\cal D}=i({\cal A})$ is {\bf involutive}. Then, $({\cal A},i, M,[.,.,]_{|{\cal S}({\cal A})})$ is a Banach Lie algebroid. So it follows from Theorem 5 that ${\cal D}$ is an integrable distribution. Thus, we get  a version of Frobenius Theorem,  as we can find in \cite{Gl}, when  $\Pi: {\cal A} \ap M$ is a  (closed) subbundle of $ TM\ap M$. In  the previous  general situation, we  can also  consider   this result as an appropriate version of {\bf  Frobenius Theorem}.
\item  Let  $\Pi: {\cal A} \ap M$ be a Banach bundle  and  $\Psi: {\cal A} \ap TM$   an injective morphism bundle. If
${\cal D}=\Psi({\cal A})$ satisfies the condition (LB) of Theorem 4, then ${\cal D}$ is integrable. From the injectivity of $\Psi$, it follows that we can define a Lie  algebra structure on the  sections ${\cal S}({\cal A})$, by:
$$\{s_1,s_2\}=\Psi^{-1}([\Psi(s_1),\Psi(s_2)]$$
So, we get a  Banach Lie algebroid  structure on $\cal A$.
\item Consider a smooth right action $\psi:M \times G\ap M$ of a connected  Banach Lie group $G$ over a Banach manifold $M$. Denote  by $\cal G$ the Lie algebra of $G$. We have a natural morphism $\xi$ of Lie algebra from $\cal G$ to ${\cal X}(M)$ which is defined  by
$$\xi_X(x)=T_{(x,e)}\psi(0,X)$$
For any $X$ and $Y$ in $\cal G$, we have:
$$\xi_{\{X,Y\}}=[\xi_X,\xi_Y]$$
 where $\{\;,\;\}$ denote the Lie algebra bracket on $\cal G$ (see for instance \cite{KrMi} chap. VIII, 36.12 or \cite{Bo}).
 On the trivial bundle $M\times {\cal G}$, each section can be identified with a map $\s: M\ap {\cal G}$  we define a Lie bracket on the set of sections by
 $$\{\{\s,\s'\}\}(x)=\{\s(x),\s'(x)\}+d\s(\xi_{\s'(x)})-d\s'(\xi_{\s(x)})$$
According to the triviality of $M\times {\cal G}$, we get a localizable Lie bracket. An anchor morphism $\Psi: M\times {\cal G}\ap TM$ is defined  by  setting $\Psi(x,X)=\xi_X(x)$.\\
  It follows that  $(M\times {\cal G},\Psi,M,\{\{\;,\;\}\})$ has  a  Banach  Lie algebroid structure  on $M$.
 \\
 Denote by $G_x$ the closed subgroup of isotropy of a point $x\in M$ and ${\cal G}_x\subset {\cal G}$ its Lie sub-algebra. Of course, we have $\ker\Psi_x={\cal G}_x$. According to Theorem 4,  if ${\cal G}_x$ is complemented in ${\cal G}$ for any $x\in M$, the weak distribution ${\cal D}=\Psi(M \times  {\cal G})$ is integrable. In fact the  leaf through $x$ is its orbit $\psi(x,G)$.
\end{enumerate}
  \end{exe}

\noindent\begin{proof}\so{\it Proof of Theorem 5}${}$\\

We will show that the property (LB) of Theorem 4 is satisfied in our context. Of course  the set ${\cal S}({\cal A})$ of local  sections of $\cal A$ is a generating upper set. Suppose that we have a structure of Banach Lie algebroid  on ${\cal A}$.  As (LB) is a local property, we may assume that $M$ is an open set of $E$ and ${\cal A}\equiv F\times M$ if $F$ is the typical fiber of ${\cal A}$. So we adopt the (local) notation used in the proof of Theorem 4.

Consider any  section $\s\in {\cal S}({\cal A})$   and fix some $x\in V$. Again, we set $Z_\s=\Psi\circ \s$ and $Z_u=\Psi(u, )$ an upper section. Given an integral curve $\g(t)=\phi^{Z_\s}_t(x)$ on $]-\varepsilon,\varepsilon[$, from (\ref{alge}), we have
$$[Z_\s,Z_u](\g(t))=\Psi(\{\s, s_u\}(\g(t))) \textrm{ where } s_u(x)=(u,x).$$
 But, using the same arguments as the ones  used in the proof of Lemma \ref{compLie}, we can show that the map
$$t\mapsto \{\s, s_u\}(\g(t))$$ is a smooth field of endomorphisms of $F$. It follows that ${\cal D}$ satisfies (LB), and then, ${\cal D}$ is integrable.\\
\end{proof}

%%%%%%%%%%%%%%%%%%%%%%%%%%%%%%%%%%%%%%%%%
\subsection{ Banach Poisson manifold}\label{balipo}
%%%%%%%%%%%%%%%%%%%%%%%%%%
We first recall the context of  Banach Poisson manifold  studied these last years (see for example \cite{OdRa2}). In particular, we will prove in a large context the existence of weak symplectic leaves.\\

A {\bf  Lie bracket } on $C^\infty(M)$ is $\R$-bilinear antisymmetric pairing $\{\;,\;\}$ on $C^\infty(M)$ which satisfies the Leibniz rule:  $\{fg,h\}=f\{g,h\}+g\{f,h\}$   and the Jacobi identity:\\ $\{\{f, g\}, h\} + \{\{g, h\}, f\} + \{\{h, f\}, g\} =0$ for  any $f,g,h\in C^\infty(M)$.

A {\bf Poisson anchor}  on $M$ is a bundle morphism $\Psi:T^*M\ap TM$  which is antisymmetric (i.e. such that  $<\a,\Psi\b>=-<\b,\Psi\a>$ for any $\a,\b\in T^*M$).

We can associate to a such morphism an $\R$-bilinear antisymmetric pairing $\{\;,\;\}$ on the set ${\cal A}^1(M)$ of $1$-form  on $M$ defined by:
$$\{\a,\b\}_\Psi=<\b,\Psi\a>$$
Moreover, for any $f \in C^\infty(M)$ we have:
$$\{f\a,\b\}_\Psi= <\b,f\Psi\a>=f\{\a,\b\}$$

\noindent So,  we get a bracket $\{\;,\;\}_\Psi$ on $C^\infty(M)$ defined by
$$\{f,g\}=\{df,dg\}_\Psi$$
As in finite dimension, $\{\;,\;\}_\Psi$ satisfies the  Jacobi identity if and only if the  the  Schouten-Nijenhuis bracket $[P,P]$ of $P$ is identically zero (see for instance  \cite{MaRa}). In this case, $C^\infty(M)$  has a structure of Lie algebra and $(M,\{\;,\;\}_\Psi)$ is called a {\bf Banach  Poisson manifold}  (see for instance \cite{OdRa1} or  \cite{OdRa2}).\\

{\it  From now, the Poisson anchor $\Psi$ is fixed and for simplicity we denote $\{\;,\;\}$ the Lie bracket associated to $\Psi$.}\\

Given such a Banach Poisson manifold, the distribution ${\cal D}=\Psi(T^*M)$ is called the {\bf characteristic distribution}. Of course, in general ${\cal D}$ is not a closed distribution but it is a weak distribution.

\noindent Associated to $\{\;,\;\}$, on $T^*M$, we have a natural skew-symmetric bilinear form $\o$ defined as follows:\\
 for any $\a$ and $\b$ in $T_x^*M$, we have $\o(\a,\b)=\{f,g\}$ if $f$ and $g$ are smooth functions defined on a neighbourhood of $x$ and such that $df(x)=\a$ and $dg(x)=\b$ (this definition is independent of the choice of $f$ and $g$).
 \\
 For each $x$, on the quotient $T_x^*M/ \ker\Psi_x$ we get a skew-symmetric bilinear form $\hat{\o}_x$. On the other hand, let $\hat{\Psi}_x:T_x^*M/ \ker\Psi_x \ap \tilde{\cal D}_x$ be the canonical isomorphism associated to $\Psi_x$ between Banach spaces. Finally we get a skew-symmetric bilinear form $\tilde{\o}_x$ on $\tilde{\cal D}_x$ such that :
 $$[\hat{\Psi}_x]^*\tilde{\o}_x=\hat{\o}_x$$

 According to \cite{OdRa2}, a {\bf symplectic leaf} of ${\cal D}$ is a weak submanifold  $({\cal L},i)$ where ${\cal L}\subset M$ and $i:{\cal L}\ap M$ is the natural inclusion with the following properties:
\begin{enumerate}
\item[(i)] $({\cal L},i)$ is a maximal integral manifold of ${\cal D}$ (in the sense of Theorem 1 part (2));
\item[(ii)] on ${\cal L}$ we have a weak symplectic form $\o_{\cal L}$ such that $(\o_{\cal L})_x=\tilde{\o}_x$ for all $x\in {\cal L}$
\end{enumerate}

\begin{rem}${}$\\
As in the context of Lie Banach algebroids,  we will say that a Lie bracket $\{\;,\;\}$ on $C^\infty(M)$ is {\bf localizable} if $\{\;,\;\}$  is compatible with the sheaf of germs of functions on $M$ .
 From our definition of Banach Poisson manifold, the Lie bracket associated to a Poisson anchor $\Psi$ is always localizable. On the other hand, given any Lie bracket $\{\;,\;\}$ on $C^\infty(M)$,   when $M$ is regular, we can prove that $\{\;,\;\}$ is localizable, and then we have a morphism $\Psi:T^*M\ap T^{**}M$ naturally associated (see \cite{CaPe}). If moreover, $\Psi(T^*M)\subset TM$, then we get the previous definition of Banach Poisson manifold (see for instance \cite{OdRa1} or  \cite{OdRa2}).
\end{rem}

\begin{theor}\label{thVI} ${}$\\
Let be $\Psi: T^*M\ap TM$ a Poisson anchor.  If the kernel of $\Psi$ is complemented in each fiber, then the associated characteristic distribution ${\cal D}$ is integrable. Moreover, each maximal integral manifold has a natural structure of weak symplectic  leaf.
\end{theor}

\noindent For an illustration of this result, the reader will find many examples of   Banach Poisson manifolds in \cite{OdRa1} and \cite{OdRa2}.\\

\noindent\begin{proof}\so{\it Proof of Theorem 6}${}$\\

At first, we can observe that  the set 
$${\cal S} =\{\Psi(df) : f\in C^{\infty}(U), U \textrm{ any open set in } M\}$$
 is an upper generating set for ${\cal D}$: given any $x\in M$, modulo any local chart  around $x$,  we can suppose that $M$ is an open subset of $ E$ and $T^*M\equiv E^*\times M$; for any $\a\in E^*$ the function $f_\a(x)=<\a,x>$ is a smooth map on $M$ such that $df_\a(y)=\a$ for any $y\in M$; so $Z_\a=\Psi(\a,y)=\Psi(df_\a(y))$  is an upper section.\\
 For any smooth  local function $f:U\ap \R$, we set $Z_f=\Psi(df,\;)$. From the Jacobi identity in $C^\infty(M)$ we have
 \begin{eqnarray}\label{liemorph}
[Z_f,Z_g]=\Psi(d\{f,g\},\;) \textrm{ for any } f,g\in C^\infty(M)
\end{eqnarray}

 According to Theorem 4,  to prove the integrability of $\cal D$, we have only to prove (LB) for the generating upper set $\cal S$. As (LB) is a local property, again we assume that $M$ is an open set in $E$. So fix some smooth function $f:M\ap \R$ and consider an integral curve $\g(t)=\phi^{Z_f}_t(x)$ through $x\in M$ defined on $]-\varepsilon,\varepsilon[$. For any $\a\in E^*$ , using  (\ref{liemorph}), we have:
 $$[Z_f,Z_\a](\g(t))=[\Psi(df),\Psi(df_\a)](\g(t))=\Psi(d\{f,f_\a\}(\g(t)))$$
But, using the same arguments as the ones used in the proof of Lemma \ref{compLie} we can show that the map
$$y \mapsto [\a \mapsto d\{f ,f_\a\}(y)]$$
is a smooth  field of continuous operators from $E^*$ to $E^*$. It follows that ${\cal D}$ satisfies (LB), and then, ${\cal D}$ is integrable.\\

Assume now that ${\cal D}$ is integrable and choose any maximal leaf ${\cal L}$. As $T_x{\cal L}=\tilde{\cal D}_x$,  on $T_x{\cal L}$ we have  the skew-symmetric bilinear form
$\tilde{\o}_x$ previously defined. We will show that $\tilde{\o}_x$ defines a closed $2$-form $\o_{\cal L}$ on $\cal L$, which is a weak symplectic form.\\

Fix  $x\in {\cal L}$. We have $T_x^*M=\ker \Psi_x\oplus S_x$.  So ${\cal L}$ is a Banach manifold modelled on $S_x$.  From the definition of $\tilde{\o}_x$, we have
\begin{eqnarray}\label{omega}
\tilde{\o}_x (\tilde{\theta}_x(\a),\tilde{\theta}_x(\b))=<\a,\tilde{\theta}_x(\b)>.
\end{eqnarray}
As we know that $\tilde{\theta}_x$ is an isomorphism from $S_x$ to $T_x{\cal L}$   it follows that $\tilde{\o}_x$ is a weak symplectic $2$-form on the Banach space $T_x{\cal L}$. On one  hand, locally, in keeping to Lemma \ref{Xdiff}, it follows that $\o_{\cal L}$ defined by $(\o_{\cal L})_x=\tilde{\o}_x$ is a smooth  differential $2$-form on ${\cal L}$. On the other hand,  for any smooth function $f$ defined on an open set $U\subset M$, we set $\tilde{f}=f\circ i$. So for any smooth functions $f$, $g$ and $h$ defined on $U$, the Jacobi identity is satisfied for $\tilde{f}$, $\tilde{g}$ and  $\tilde{h}$ on the open set  $i^{-1}(U)\subset {\cal L}$. So, by classical arguments of Poisson bracket (see for instance \cite{LiMa}, \cite{OdRa1}, \cite{OdRa2}), we get:
$d\o_{\cal L}(i_*Z_f,i_*Z_g,i_*Z_h)=0$ for any choice of functions $f$, $g$ and  $h$. So $\o_{\cal L}$ is closed and the proof of Theorem 6 is complete.

\end{proof}

%\begin{thebibliography}{99}

\end{document}